\documentclass[10pt,leqno]{amsart}
\usepackage{graphicx}
\usepackage{indentfirst,csquotes,comment}

\usepackage{multicol}
\newcommand{\AuthorBlock}[3]{%
  \noindent\textbf{#1}\\%
   #2\\%
  E-mail: \texttt{#3}%
  \par\vspace{1em}%
}
\makeatletter
\patchcmd{\@settitle}
  {\uppercasenonmath\@title}
  {}
  {}
  {}

\patchcmd{\@setauthors}
  {\MakeUppercase{\authors}}
  {\authors}
  {}
  {}
\makeatother

\usepackage{amssymb,amsthm,amsmath,amsfonts}
\usepackage{xcolor,paralist,hyperref,fancyhdr,etoolbox,lipsum}
\newtheorem{theorem}{Theorem}[section]
\newtheorem{definition}[theorem]{Definition}
\newtheorem*{introdefinition}{Definition}

\newtheorem{lemma}[theorem]{Lemma}
\newtheorem{proposition}[theorem]{Proposition}

\newtheorem{conjecture}[theorem]{Conjecture}
\newtheorem*{remark}{Remark}

\newcommand{\D}{\mathbb{D}}  
\newcommand{\Z}{\mathbb{Z}}
\newcommand{\C}{\mathbb{C}}  
\newcommand{\N}{\mathbb{N}}   
\newcommand{\Q}{\mathbb{Q}}    
\newcommand{\PP}{\mathbb{P}}  
\newcommand{\EE}{\mathbb{E}} 
\newcommand{\E}{\mathbb{E}}   

\newcommand{\IIA}{\mathrm{IAG}}
\newcommand{\ICP}{\mathrm{ICP}}

\newcommand{\IIP}{\mathrm{IHP}}

\newcommand{\eps}{\varepsilon}

\newcommand{\Hh}{\mathbb{H}}           
\newcommand{\1}{\mathbf{1}}

\newcommand{\bdry}{\partial}
\newcommand{\supp}{\operatorname{supp}}
\newcommand{\Pp}{\mathbb{P}}

\hypersetup{ colorlinks=true, linkcolor=black, filecolor=black, urlcolor=black }

\begin{document}
\title[Unimodular random ideal hyperbolic polyhedra: geometric characteristics and random walks]{Unimodular random ideal hyperbolic polyhedra: geometric characteristics and random walks}
\author[Ge, Lu, Wang and Zhou]
{Huabin Ge, Yangxiang Lu, Chuwen Wang, and Tian Zhou}

\let\thefootnote\relax

\begin{abstract}
We initiate a systematic study of unimodular random infinite ideal hyperbolic polyhedra in three-dimensional hyperbolic space, developing a geometric and probabilistic theory of their global behavior. At the heart of this theory is an intrinsic geometric characteristic, determined entirely by the dihedral angles surrounding a polyhedral face. We establish a unimodular Gauss–Bonnet identity and show that the expectation of this local geometric invariant completely determines the global conformal type: under suitable assumptions, a random ideal polyhedron is almost surely parabolic precisely when the expected characteristic vanishes, and hyperbolic precisely when it is negative. This provides a sharp geometric classification of random infinite ideal polyhedra.\\[-4pt]

Building on the correspondence with ideal circle patterns, we further develop a refined Ring Lemma that yields quantitative geometric control without bounded-degree assumptions. Combining these geometric tools with unimodular probability methods, we identify the Poisson boundary of the face random walk and establish positive linear hyperbolic speed. Together, these results establish a local-to-global framework linking intrinsic dihedral geometry, conformal uniformization, and stochastic behavior at infinity for unimodular random infinite ideal hyperbolic polyhedra.

\end{abstract}

\maketitle
\bigskip
\tableofcontents

\parskip=3pt plus 2pt minus 1pt  

\section{Introduction}

The theory of hyperbolic polyhedra has a long history, beginning with the foundational work of Andreev, Thurston, and Rivin on existence, rigidity, and the realization of prescribed combinatorial and dihedral-angle data \cite{And70,Koebe36,Riv96,ThurNotes}. For infinite hyperbolic polyhedra, however, fundamentally new questions arise. Beyond existence and rigidity, one may ask how local geometric information governs the global structure of a polyhedron and its behavior at infinity. When the combinatorial structure and dihedral angles are random, these questions lead naturally to an interaction between three-dimensional hyperbolic geometry, discrete conformal geometry, and probability.

In this paper, we initiate a systematic study of \textit{unimodular random infinite ideal hyperbolic polyhedra} in hyperbolic three-space. Our central problem is to understand whether the global geometry and stochastic behavior of such a polyhedron can be characterized through geometric quantities determined locally by its three-dimensional dihedral-angle data. We establish a theory that provides an affirmative answer. At its heart is an intrinsic geometric characteristic associated with each polyhedral face. We prove that its expectation satisfies a unimodular Gauss–Bonnet identity and, under suitable geometric and integrability assumptions, completely determines the global conformal type of the random polyhedron. More precisely, the expected characteristic distinguishes two geometric regimes: the parabolic regime, in which the ideal vertices accumulate at a single point at infinity, and the hyperbolic regime, in which their accumulation set is a circle. Thus, the global geometry of an infinite random polyhedron is classified by the expectation of a quantity determined entirely by its local dihedral geometry.

Our theory extends beyond this geometric classification. In the hyperbolic regime, we study the asymptotic behavior of simple random walk on the face-adjacency graph of the polyhedron. We establish convergence to a geometric boundary, identify the Poisson boundary, and prove positive linear hyperbolic speed. A crucial ingredient is a refined Ring Lemma for ideal circle patterns with prescribed intersection angles, providing quantitative geometric control without assuming uniformly bounded degrees.

The central step in our approach is to develop a geometric mechanism that transforms the intrinsic dihedral-angle data of infinite ideal hyperbolic polyhedra into information about their global geometry and stochastic behavior. Building on the existence, rigidity, and uniformization theory of ideal circle patterns \cite{GeIIP2,GHZ21ICP,GeIIP1}, we establish a unimodular Gauss-Bonnet principle for geometric characteristics, through which local three-dimensional geometry governs the global conformal type. Equally important, we develop new quantitative geometric estimates for ideal circle patterns, including a refined Ring Lemma that controls intersecting-circle configurations without imposing bounded-degree assumptions. These geometric developments provide the essential bridge to the theory of unimodular random planar maps and the mass transport principle \cite{AldLy07,AHNR16,map}, enabling us to establish both the sharp parabolic–hyperbolic classification and the stochastic boundary theory of random infinite ideal hyperbolic polyhedra.

\subsection{Random infinite ideal hyperbolic polyhedra and geometric characteristics.}
Let $P\subset\mathbb H^3$ be an infinite \textit{ideal hyperbolic polyhedron} (IHP). We consider convex polyhedra with infinitely many faces and ideal vertices on the sphere at infinity $\partial_\infty\mathbb H^3$, in the locally finite setting specified in Section \ref{section-iip-H3}. Unlike a finite ideal polyhedron, whose ideal vertices form a finite subset of the sphere at infinity, an infinite ideal polyhedron necessarily has accumulation points. The geometric structure of these accumulation points provides a natural way to distinguish different types of infinite polyhedra. The existence and rigidity theory of infinite ideal circle patterns developed in \cite{GeIIP2,GHZ21ICP,GeIIP1} leads to two fundamental types: an infinite ideal polyhedron is called \emph{parabolic} if its ideal vertices have exactly one accumulation point on $\partial_\infty\mathbb H^3$, and \emph{hyperbolic} if their accumulation set is a circle and all ideal vertices lie on the same side of this circle. These notions correspond to the parabolic and hyperbolic uniformization types of the associated ideal circle patterns; precise definitions are given in Section \ref{section-iip-H3}.

We study random rooted ideal polyhedra $(P,f_\rho)$, where $f_\rho$ is a distinguished face. Unimodularity is defined through the rooted marked dual graph of $P$, whose vertices represent polyhedral faces and whose edge marks record the corresponding dihedral angles. This allows the mass transport principle to be applied while retaining the essential three-dimensional geometric information. The central quantity in our theory is defined as follows.

\begin{introdefinition}[geometric characteristic]
Let $P$ be an infinite ideal hyperbolic polyhedron. For each face $f$ of $F(P)$ (see Figure \ref{geom-charact} for example), define
\begin{equation}
\label{eq:ihp-characteristic}
T_P(f):=2\pi-\sum_{e\subset\partial f}\Theta_P(e),
\end{equation}
where $\Theta_{P}(e)\in(0,\pi)$ is the dihedral angle of $P$ at $e$. We call $T_P(f)$ the \emph{geometric characteristic number} of $f$ , and the collection $T_P=(T_P(f):f\in F(P))$ the \emph{geometric characteristic} of $P$.
\end{introdefinition}

\begin{figure}[h]
\centering
\includegraphics[scale=0.18]{fig/1027.ai}
   \caption{The geometric characteristic under the IHP/IAG duality.
    Left: a face $f$ of the IHP ($T_P(f)=2\pi-\sum_{k=1}^{4}\Theta_P(e_k)$).
    Right: the corresponding vertex $i=f^*$ of the IAG $(G,\Theta)$.
($T(i)=2\pi-\sum_{k=1}^{4}\Theta(e_k^*)=T_P(f)$).}
\label{geom-charact}
\end{figure} 

The geometric characteristic is intrinsic to the dihedral geometry of the polyhedron: it is determined by the angles between its supporting hyperbolic planes and is independent of the choice of circle-pattern radii or auxiliary coordinates. This notion is closely related to the geometric-characteristic theory of circle patterns developed in \cite{GeIIP2,GeLin2024,LLS25}. Its importance here lies in the possibility that a simple local geometric quantity may encode global information about an infinite polyhedron. To see its motivation, consider first a finite convex ideal hyperbolic polyhedron. The ideal-vertex angle relations, together with Euler's formula, imply
\begin{equation}
\label{eq:finite-gauss-bonnet}
\sum_{f\in F(P)}T_P(f)=4\pi.
\end{equation}
This is a Gauss--Bonnet identity expressed entirely in terms of the dihedral geometry of the polyhedron. For an infinite polyhedron, however, the corresponding global sum need not be meaningful, and a local geometric characteristic does not generally determine the conformal type of an arbitrary deterministic realization. Unimodularity provides a natural substitute for global summation: expectation at a random root.

To formulate the resulting identity, we recall
the combinatorial curvature used in the theory
of unimodular random planar maps \cite{map}:
\begin{equation}\label{eq:ihp-combinatorial-curvature}
\kappa_P(f):=2\pi-\sum_{v\in\partial f}
\frac{(\deg_P(v)-2)\pi}{\deg_P(v)},
\end{equation}
where $\deg_P(v)$ denotes the number of edges incident to the ideal vertex $v$. Our first main result establishes an infinite-volume Gauss-Bonnet identity for random ideal polyhedra.

\begin{theorem}[Unimodular Gauss-Bonnet identity]
\label{thm:ihp-gb}
Let $(P,f_\rho)$ be a unimodular random rooted infinite ideal hyperbolic polyhedron. Then
\begin{equation}
\label{eq:ihp-gb}
\mathbb E[T_P(f_\rho)]=\mathbb E[\kappa_P(f_\rho)].
\end{equation}
Moreover, both expectations are finite if $\mathbb E[\deg(f_\rho)]<\infty$.
\end{theorem}

This identity reveals a remarkable interaction between geometry and probability. The characteristic number $T_P(f)$ depends on the actual dihedral angles of the polyhedron, whereas $\kappa_P(f)$ depends only on its combinatorial structure. Under unimodularity these two quantities agree in expectation: the mass transport principle converts the angle relations at individual ideal vertices into an infinite-volume curvature identity. In particular, explicit dependence on the angles disappears after averaging. This is striking in view of the deterministic theory of infinite ideal circle patterns, where the conformal type may depend on the prescribed intersection angles. Nevertheless, the averaged characteristic retains decisive geometric information, as our principal classification theorem shows.

\begin{theorem}[Geometric dichotomy]
\label{thm:ihp-dichotomy}
Let $(P,f_\rho)$ be an infinite, tame, ergodic, unimodular random rooted ideal hyperbolic polyhedron. Suppose that $\sup_e\Theta_P(e)<\pi$ almost surely, and assume $\mathbb E[\deg(f_\rho)]<\infty$. Then $\mathbb E[T_P(f_\rho)]\leq 0$ and 
\begin{equation}\label{eq:ihp-dichotomy}
\begin{aligned}
\mathbb E[T_P(f_\rho)]=0
&\quad\Longleftrightarrow\quad P\text{ is almost surely parabolic},\\
\mathbb E[T_P(f_\rho)]<0
&\quad\Longleftrightarrow\quad P\text{ is almost surely hyperbolic}.
\end{aligned}
\end{equation}
\end{theorem}

The geometric dichotomy theorem provides a sharp geometric classification of unimodular random infinite ideal hyperbolic polyhedra. A global distinction concerning the accumulation of infinitely many ideal vertices on the sphere at infinity is characterized by the expected value of a local quantity computed from three-dimensional dihedral angles. Its significance goes beyond a reformulation of combinatorial curvature: it identifies a characteristic naturally associated with the polyhedron itself and shows that, after unimodular averaging, this characteristic governs the global uniformization type. Together, the unimodular Gauss–Bonnet identity and the geometric dichotomy establish a local-to-global principle for random infinite polyhedral geometry.

\subsection{Ideal circle patterns and the geometry of polyhedral realization}
\begin{figure}[h]
\centering
\includegraphics[scale=0.018]{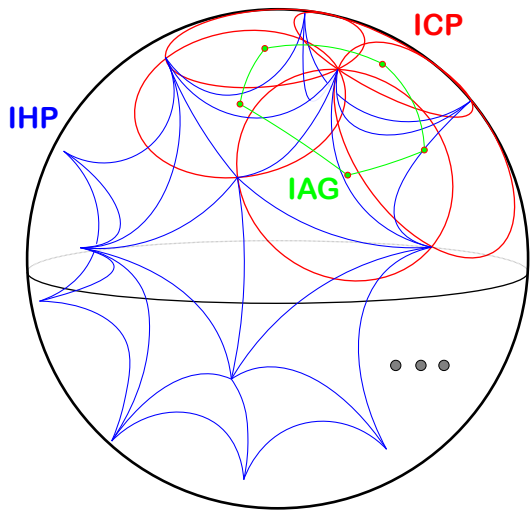}
\caption{$\ICP$-$\IIP$ correspondence and \(\IIP/\IIA\) duality.}
\label{ICP-IHP}
\end{figure}

The connection between three-dimensional hyperbolic polyhedra and planar geometry is provided by circles on the sphere at infinity \cite[Chapter 13]{ThurNotes}. In the Poincar\'e ball model, the hyperbolic plane containing each face of an IHP $P$ determines a circle on $\partial_\infty\mathbb H^3$. When two faces meet along an edge, the corresponding circles intersect at an angle equal to the dihedral angle of the polyhedron; the circles associated with faces incident to an ideal vertex pass through a common point at infinity. The resulting configuration is an \emph{ideal circle pattern} (ICP). Conversely, under suitable admissibility conditions, an ICP determines an IHP through the intersection of the associated hyperbolic half-spaces. This is the ICP-IHP correspondence illustrated in Figure \ref{ICP-IHP}.

A finite prototype is Rivin's characterization of convex ideal hyperbolic polyhedra \cite{Riv96}. Let $\mathcal{D}=(V,E,F)$ be a finite cellular decomposition of the sphere, equipped with an angle assignment $\Theta:E\to(0,\pi)$. Rivin's theorem states that there exists a convex ideal hyperbolic polyhedron with dual combinatorial type $\mathcal{D}$ and prescribed dihedral angles $\Theta$ if and only if
\begin{align}
\tag{$C_1$}\label{eq:C1}
\sum_{e\in\partial f}(\pi-\Theta(e))&=2\pi,
&&\text{for every face $f\in F$, and}\\
\tag{$C_2$}\label{eq:C2}
\sum_{e\in\gamma}(\pi-\Theta(e))&>2\pi,
&&\text{for every cycle $\gamma$ not bounding a single face}.
\end{align}
Moreover, the polyhedron is unique up to isometry. Condition $(C_1)$ expresses the angle relation imposed by an ideal vertex, while $(C_2)$ is a global convexity condition. For infinite cellular decompositions, we use the corresponding Rivin-type conditions and the strengthened uniform cycle-separation condition
\begin{equation}
\tag{$C_2'$}
\label{eq:C2prime}
\inf_{\gamma\neq \partial f}\Big\{\sum_{e\in\gamma}(\pi-\Theta(e))\Big\}>2\pi,
\end{equation}
where the infimum is taken for all cycle $\gamma$ which does not bound a single face. This strengthened condition is useful for the quantitative estimates in our boundary theory.

\begin{introdefinition}[IAG]\label{def:iag}
Let $\mathcal{D}=(V,E,F)$ be an infinite, locally finite disk cellular decomposition, and let $G=(V,E)$ be its underlying graph. An angle assignment $\Theta:E(G)\to(0,\pi)$ defines an \textit{ideal angled graph} (IAG) $(G,\Theta)$ if $(C_1)$ and $(C_2)$ hold. If $(C_2)$ is strengthened to $(C_2')$, the IAG is called \emph{tame}. A rooted IAG is a triple $(G,\rho,\Theta)$, where $\rho\in V(G)$ is the distinguished root.
\end{introdefinition}

Throughout the paper, IAGs are assumed to be simple, infinite, locally finite, and one-ended. The existence and rigidity theory developed by Ge-Hua-Zhou and Ge-Yu-Zhou \cite{GeIIP2,GHZ21ICP,GeIIP1} provides the corresponding infinite ICP realization. Under the appropriate angle hypotheses, an IAG admits an ICP realization in either $\mathbb C$ or $\mathbb D$, with the corresponding rigidity properties. Through the ICP-IHP correspondence, this realization determines an infinite IHP $P(G,\Theta)\subset\mathbb H^3$ whose dual $1$-skeleton is $G$ and whose dihedral angles are prescribed by $\Theta$. Conversely, an infinite IHP in our setting determines an IAG through its face-adjacency structure and dihedral-angle assignment. We refer to this relationship as the \emph{IHP/IAG duality} (Figure \ref{geom-charact}).

This duality makes unimodular probability methods available for the dual planar structure while preserving the interpretation of the results in genuine three-dimensional polyhedral geometry. Its prescribed angles are not arbitrary weights: they satisfy the compatibility relations imposed by ideal vertices, and these relations are essential to existence, rigidity, and quantitative control. Thus the ideal circle patterns used here are geometrically constrained realizations, rather than ordinary tangency packings with freely introduced angle parameters.

\subsection{Unimodular ideal angled graphs and the classification theorem.}
\label{subsec:iag-dichotomy}
We now turn to the probabilistic formulation of the geometric dichotomy. The theory of unimodular random rooted graphs developed by Aldous-Lyons \cite{AldLy07} provides a framework for infinite random graphs and their local weak limits, based on the mass transport principle. For unimodular planar maps, Angel-Hutchcroft-Nachmias-Ray \cite{AHNR16,map} established fundamental relationships among combinatorial curvature, invariant amenability, circle-packing type, vertex extremal length, and random walks. Their work provides important probabilistic tools for the classification below. In our setting the dual structure is an angle-marked planar map whose marks are actual polyhedral dihedral angles; the natural local quantity is therefore
\begin{equation}\label{eq:iag-characteristic}
T(v):=2\pi-\sum_{e\ni v}\Theta(e),
\end{equation}
which is related to the angle-sum character
studied in \cite{LLS25},
with $T(v)=2\pi-\mathcal L(v)$;
see also \cite{GeIIP2}. Under the IHP/IAG duality, $T(v)$ is precisely the characteristic number of the corresponding polyhedral face (see Figure \ref{geom-charact}). For comparison, the combinatorial curvature at $v$ is
\begin{equation}\label{eq:iag-curvature}
\kappa(v):=2\pi-\sum_{f\ni v}
\frac{(\deg(f)-2)\pi}{\deg(f)}.
\end{equation}

Theorem~\ref{thm:ihp-gb} has the equivalent graph-theoretic form
\begin{equation}\label{eq:iag-gb}
\mathbb E[T(\rho)]=\mathbb E[\kappa(\rho)],
\end{equation}
connecting intrinsic ideal-polyhedral geometry with curvature theory for unimodular planar maps.

\begin{theorem}[Dichotomy for unimodular IAGs]
\label{thm:iag-dichotomy}
Let $(G,\rho,\Theta)$ be an ergodic unimodular random tame rooted IAG with $\mathbb E[\deg(\rho)]<\infty$. Then $\mathbb E[T(\rho)]\leq0$, and the following three conditions are equivalent, with the parenthetical alternatives corresponding to each other:
\begin{enumerate}
\item $\mathbb E[T(\rho)]=0$ (respectively, $\mathbb E[T(\rho)]<0$);
\item $G$ is invariantly amenable (respectively, invariantly non-amenable);
\item $G$ is almost surely VEL-parabolic (respectively, VEL-hyperbolic).
\end{enumerate}
If in addition $\sup_{e}\Theta(e)<\pi$, these conditions are equivalent to
\begin{enumerate}
\setcounter{enumi}{3}
\item $G$ is almost surely ICP-parabolic (respectively, ICP-hyperbolic).
\end{enumerate}
If $G$ further has bounded degree, they are also equivalent to
\begin{enumerate}
\setcounter{enumi}{4}
\item simple random walk on $G$ being recurrent (respectively, transient).
\end{enumerate}
Moreover, almost sure ICP-hyperbolicity implies $\mathbb E[T(\rho)]<0$ even without a bounded-degree assumption or the uniform angle bound.
\end{theorem}

The theorem identifies several aspects of the  random geometric structure. Invariant amenability concerns its probability law, VEL captures potential-theoretic type, ICP uniformization describes conformal geometry, and the geometric characteristic is read directly from the dihedral angles of the polyhedron. In addition to the above equivalences, we also have a list of 17 equivalences as in Theorem 1 in \cite{map}, relating to aspects of the map including amenability, random walks, harmonic functions, spanning forests, Bernoulli bond percolation, and the conformal type of associated Riemann surfaces. A special case occurs when $\sup_e\Theta(e)\leq\pi/3$. Condition $(C_1)$ forces the cellular decomposition to be a triangulation with $\Theta(e)=\pi/3$; consequently,
$$T(\rho)=\frac{\pi}{3}(6-\deg(\rho)),$$
and
$$\mathbb E[T(\rho)]=\frac{\pi}{3}\bigl(6-\mathbb E[\deg(\rho)]\bigr).$$
The theorem then recovers the degree-six dichotomy for unimodular random triangulations \cite{AHNR16}, with a geometric interpretation in terms of ideal hyperbolic polyhedra. For general IAGs, by contrast, $T(\rho)$ depends on the prescribed dihedral angles and is not determined by degree alone; nevertheless, its expectation obeys the universal identity \eqref{eq:iag-gb}. This interplay between local geometric dependence and unimodular averaging is a central feature of the theory.

  \subsection{Random walks, geometric boundaries, and a refined Ring Lemma}\label{subsec:boundary-introduction}

The geometric characteristic determines the global conformal regime. In the hyperbolic case, we further ask how random walks on the face-adjacency graph approach infinity. Let $(G,\rho,\Theta)$ be an ICP-hyperbolic unimodular random IAG, and let $\mathcal P$ be its ICP realization in the unit disk $\mathbb D$. We consider simple random walk $(X_n)_{n\geq0}$ on $G$, started at $\rho$; under the IHP/IAG duality, this is the face random walk on the associated polyhedron. Let $z(v)$ and $z_h(v)$ denote the Euclidean and hyperbolic centers of the circle corresponding to $v$. We study whether these centers converge to $\partial\mathbb D$ along the walk, whether this boundary realizes the Poisson boundary, and whether the walk escapes at positive hyperbolic speed.

The main obstacle in establishing these boundary results is geometric. Unlike ordinary circle packings, ideal circle patterns involve prescribed nonzero intersection angles and are subject to the compatibility
conditions imposed by ideal vertices of hyperbolic polyhedra. Moreover, the underlying random graphs may have unbounded degrees. Consequently, the
existing qualitative existence and rigidity theory does not provide the quantitative geometric control needed for random-walk analysis. A key challenge is therefore to obtain effective estimates for neighboring
circle radii that remain valid without bounded-degree assumptions. Our main geometric input is the following refined Ring Lemma, which replaces a global degree bound by local combinatorial complexity.

\begin{theorem}[Refined Ring Lemma for ICP]\label{thm:icp-refined-ring}
Let $(G,\Theta)$ be a tame IAG with $\sup_e\Theta(e)<\pi$, and let $r:V(G)\to(0,\infty)$ denote the radius function of an ICP realization. There exists a constant $C=C(G,\Theta)>0$, such that for every edge $u\sim v$,
\begin{equation}
\label{eq:icp-refined-ring}
\frac{r(v)}{r(u)}>\exp\bigl(-C S(u)\bigr),
\qquad S(u):=\sum_{w\sim u}\deg(w).
\end{equation}
The quantity $S(u)$ is called the \emph{flower degree} of $u$.
\end{theorem}

The estimate (\ref{eq:icp-refined-ring}) is particularly suited to unimodular random graphs: by the mass transport principle, we have
$\mathbb E[S(\rho)]=\mathbb E[\deg(\rho)^2]$ and $\mathbb E[\deg(\rho)S(\rho)]\leq\mathbb E[\deg(\rho)^3]$. A finite third-moment assumption thus gives the integrability needed to control logarithmic changes of circle radii along the walk. The refined Ring Lemma supplies the geometric input through which local angle-constrained configurations become accessible to unimodular random-walk methods.

\begin{theorem}[Poisson boundary]
\label{thm:icp-poisson}
Let $(G,\rho,\Theta)$ be a unimodular, ICP-hyperbolic, random rooted tame IAG. Assume $\sup_e\Theta(e)<\pi$ and $\mathbb E[\deg(\rho)^3]<\infty$. Let $(X_n)$ be simple random walk started at $\rho$. Then, conditional
on $(G,\rho,\Theta)$, almost surely:
\begin{enumerate}
\item The Euclidean centers $z(X_n)$ and hyperbolic centers $z_h(X_n)$ converge to a common random point $\xi\in\partial\mathbb D$.
\item The exit distribution of $\xi$ is non-atomic and has full support on $\partial\mathbb D$.
\item The circle $\partial\mathbb D$ realizes the Poisson boundary of $(G,\Theta)$. In particular, for every bounded harmonic function $h$ on $G$, there is a bounded measurable function $g:\partial\mathbb D\to\mathbb R$ such that $h(v)=\mathbb E_v[g(\xi)]$ for all vertices $v$.
\end{enumerate}
\end{theorem}

\begin{theorem}[Positive hyperbolic speed]
\label{thm:icp-speed}
Under the assumptions of Theorem~\ref{thm:icp-poisson},
\begin{equation}\label{eq:positive-hyp-speed}
\lim_{n\to\infty}\frac{d_{\mathbb H}(z_h(\rho),z_h(X_n))}{n}
=\lim_{n\to\infty}\frac{-\log r(X_n)}{n}>0
\end{equation}
almost surely. If $(G,\rho,\Theta)$ is ergodic, the limiting speed is almost surely constant.
\end{theorem}

Through the IHP/IAG duality, these theorems describe the asymptotic behavior of the face random walk on an infinite ideal hyperbolic polyhedron. The circle $\partial\mathbb D$, arising from the ICP-hyperbolic uniformization, is a geometric realization of its Poisson boundary. Together with the characteristic dichotomy, the boundary and speed results provide a systematic description of both the global conformal regime and the stochastic behavior at infinity of unimodular random infinite ideal hyperbolic polyhedra.

\subsection{Background and related work}\label{subsec:related-work}

Our results lie at the intersection of unimodular random planar geometry, discrete conformal geometry, and three-dimensional hyperbolic polyhedra. On the geometric side, the classical characterization of convex ideal hyperbolic polyhedra is due to Rivin \cite{Riv93,Riv96}, building on the work of Andreev and Thurston. Further developments relating polyhedra, circle patterns, and variational principles include Bobenko-Springborn \cite{BobSpr04}, Ge-Hua-Zhou \cite{GHZ19}, Springborn \cite{Spring20} and Schlenker \cite{SchlHypCirc05,SchlCircSing08,SchlHypPoly}. The Koebe-Andreev-Thurston theorem \cite{Koebe36,ThurNotes} and the work of He \cite{He99}, He-Schramm \cite{HS93,HeSc95} establish fundamental connections between planar combinatorics and conformal geometry. Classical convergence results for circle packings include the
work of Rodin-Sullivan \cite{RodSul87} and
He-Schramm \cite{He96Riem}; see also Stephenson's
monograph \cite{Steph05}. The existence, rigidity, and uniformization results for infinite ICP developed in \cite{GeIIP2,GHZ21ICP,GeIIP1} provide the geometric foundation for the present study.
 
On the probabilistic side, unimodular random rooted graphs were developed by Aldous-Lyons \cite{AldLy07}; see also Lyons-Peres \cite{LyPer16}. Angel-Hutchcroft-Nachmias-Ray \cite{AHNR16,map} established relationships among combinatorial curvature, invariant amenability, conformal type, and probabilistic properties of unimodular planar maps. Their results provide important tools for our classification theorem. The present problem nevertheless involves structures arising from actual hyperbolic polyhedra: prescribed intersection angles must satisfy ideal-vertex compatibility relations, the corresponding circle configurations must admit geometric realizations, and the relevant characteristic is read from the three-dimensional dihedral angles. The unimodular Gauss--Bonnet identity connects this intrinsic quantity to planar curvature theory, while ICP uniformization translates the classification into a geometric dichotomy for infinite polyhedra.

For random-walk boundaries, the circle-packing approach of Benjamini-Schramm \cite{BS96a} and the potential-theoretic classification of He-Schramm  \cite{HS93,HeSc95} establish
fundamental connections between conformal geometry
and random walks. Angel-Barlow-Gurel-Gurevich-Nachmias \cite{ABGN16} identified the circle-packing
boundary with the Poisson boundary under suitable
assumptions, while Hutchcroft-Peres \cite{HP17} developed a unified approach through planar
embeddings, circle packings, and square tilings.
Related results based on square tilings were obtained
by Benjamini and Schramm \cite{BS96sq} and
Georgakopoulos \cite{Geo16}.
For unimodular random triangulations, related
results were established by
Angel-Hutchcroft-Nachmias-Ray \cite{AHNR16}. Other developments in stochastic hyperbolic
geometry include Curien's planar stochastic
hyperbolic triangulations ~\cite{N15} and the study of anchored expansion and
random-walk speed for Poisson-Voronoi
tessellations by
Benjamini-Paquette-Pfeffer~\cite{IEJ14}. Our boundary theory draws on these probabilistic
developments while introducing the quantitative
geometry of ideal circle patterns needed to
study random walks on infinite ideal hyperbolic
polyhedra.

In a broader context, this work is related to random infinite planar maps and their scaling limits, including the UIPT \cite{AS03}, the UIPQ \cite{CurMenMir13}, and the Brownian plane \cite{CLG14,CLG16} and the connections between random planar maps and Liouville quantum gravity \cite{GHS19}. The study of random planar maps is also informed
by the recurrence results for planar graph limits
of Benjamini-Schramm \cite{BS01} and
Gurel-Gurevich-Nachmias \cite{GGN13},
as well as Benjamini-Curien's analysis of
simple random walk on the uniform infinite
planar quadrangulation \cite{BC13}. Quantitative connections between discrete random maps
and continuum conformal geometry include convergence
under the Cardy embedding due to
Holden-Sun \cite{HS23}, results on the fractal dimension
of Liouville quantum gravity by
Ding-Gwynne \cite{DingGwynne-LQGdim-CMP2020},
and bijective approaches to percolation on random
triangulations by
Bernardi-Holden-Sun \cite{BernardiHoldenSun-Percolation-Memoirs2023};
see also the work of Gwynne-Miller on spectral dimension,
resistance, and displacement of random walks on random
planar maps \cite{GM21}. These developments illustrate the close relationship
between probabilistic properties of random maps and their
conformal geometry. The present work pursues a
complementary direction: rather than studying scaling
limits or continuum embeddings alone, we investigate
random infinite objects with genuine three-dimensional
hyperbolic realizations and determine how their intrinsic
dihedral-angle geometry governs global classification
and stochastic behavior at infinity.

In a companion paper, \emph{Geometry of Unimodular Random Hyperbolic Polyhedra} \cite{GLWZ2}, we develop related geometric-characteristic and stochastic boundary theories for trivalent hyperbolic polyhedra admitting ordinary, ideal, and hyperideal vertices. Together, these works contribute to a broader program connecting intrinsic polyhedral geometry, unimodular Gauss-Bonnet theory, and the probabilistic behavior of infinite hyperbolic structures. 

\noindent\textbf{Organization of the paper.}
Section 2 collects the geometric, combinatorial, and
probabilistic preliminaries, including ideal circle patterns,
infinite ideal hyperbolic polyhedra, unimodular random
rooted graphs, vertex extremal length, and invariant
amenability. Section 3 proves the unimodular Gauss--Bonnet
identity and the parabolic--hyperbolic dichotomy in the
language of ideal angled graphs, yielding the corresponding
classification of random infinite ideal polyhedra and its
triangulation special case. Section 4 develops the refined
Ring Lemma for ideal circle patterns and proves the
convergence and Poisson-boundary results for random walks.
Section ~5 establishes positive linear hyperbolic speed and
relates the escape rate to the decay of circle radii.
Finally, Section 6 discusses further related topics.

\medskip

\noindent\textbf{Acknowledgements.}
The authors would like to thank Puchun Zhou, Longsong Jia,
Professor Timothy Budd, and Professor Nicolas Curien
for helpful discussions. The first author expresses his sincere thanks to Professor Xin Sun and Professor Jian Ding for their valuable suggestions and comments. The second and third authors are also grateful to Professor Xin Sun for many beneficial conversations. The fourth author is especially grateful to Professor Gang Tian for his constant support and inspiration. The first author is supported by NSFC, no.12341102, no.12122119, and no.12525103.

\section{Preliminaries}

In this section we collect the geometric, combinatorial, and probabilistic
ingredients required for this paper. 

\subsection{Ideal circle pattern and ideal angled graph}
\label{section-icp}
\begin{definition}[Ideal Circle pattern \cite{GeIIP1}]
	Let $\mathcal{D} = (V, E, F)$ be an infinite cellular decomposition of the plane. 
	A circle pattern $\mathcal{P}$ associated with $\mathcal{D}$ is a collection of circles 
	$\{\mathcal{P}(v)\}_{v \in V}$ in $\mathbb{R}^2$. 
	We say that $\mathcal{P}$ is an ideal circle pattern (abbr. $\ICP$) of $\mathcal{D}$ 
	if the following conditions hold:
	\begin{enumerate}[(1)]
		\item $\mathcal{P}(u)\cap\mathcal{P}(v)\neq\emptyset$ whenever $u\sim v$.
		\item $\cap_{v<f}\mathcal{P}(v)=\{v_f\}$ for each face $f\in F$, where $v_f$ can be regarded as a dual vertex corresponding to the face $f$ in the dual graph of $\mathcal{D}$.
	\end{enumerate}
	
	We denote by $D(v)$ the closed disk bounded by $\mathcal{P}(v)$. 
	For each edge $e = \{u, v\}$, the intersection angle of the circles 
	$\mathcal{P}(u)$ and $\mathcal{P}(v)$ is denoted by $\Theta(e) = \Theta(u, v)$.
\end{definition}

\begin{figure}[htbp]
\centering
\includegraphics[scale=0.03]{fig/icpvsiag.ai}
\caption{Ideal circle pattern}
\label{4}
\end{figure}

The angle $\Theta:E\to(0,\pi)$ may satisfy the
Rivin-Thurston combinatorial conditions:

\begin{align*}
	\text{($C_1$)}&\quad
	\sum_{e\in\partial f} (\pi-\Theta(e)) = 2\pi,
	&&\forall f\in F;\\[2mm]
	\text{($C_2$)}&\quad
	\sum_{e\in\gamma} (\pi-\Theta(e)) > 2\pi,
	&&\forall \text{ cycle } \gamma \text{ not bounding a face};\\[2mm]
	\text{($C_2^\prime$)}&\quad
	\inf_{\gamma\neq\partial f}\Big\{\sum_{e\in\gamma}(\pi-\Theta(e))\Big\}>2\pi.&& \gamma  \text{ runs over all cycles not bounding a single face}.
\end{align*}

\textbf{
Throughout the paper,
whenever $(C_2')$ holds, we choose $\varepsilon_0>0$
such that
$\sum_{e\in\gamma}(\pi-\Theta(e))>2\pi+\varepsilon_0$
for every cycle $\gamma$ not bounding a face.
Whenever $\sup_{e\in E}\Theta(e)<\pi$, we choose
$\varepsilon_1>0$ such that
$\Theta(e)\le\pi-\varepsilon_1$ for every $e\in E$.
}

We present a geometric approach for constructing ideal circle patterns. 
Given a radius assignment \(r \in \mathbb{R}_+^V\), known in the literature as a \textbf{circle pattern metric}, 
we augment the cellular decomposition \(\mathcal{D}\) as follows. Let $G=(V,E)$ be its underlying graph.
For each face \(f\), we place an auxiliary vertex \(v_f\) in its interior.
The collection of these dual vertices is denoted by \(V_F\). 
The \textbf{incidence graph} of \(G\) (cf.~\cite{Coxeter1950}) is then introduced.

\begin{definition}
The incidence graph \(I(G)\) is a bipartite graph with bipartition \(\{V,V_F\}\). 
A vertex \(v \in V\) and a dual vertex \(v_f \in V_F\) are adjacent in \(I(G)\) if and only if \(v\) lies on the boundary of the face \(f\).
\end{definition}

The incidence graph demonstrates the incidence relations between vertices and faces. 
For an edge \(e=\{v,w\}\in E\), let \(f_1,f_2\in F\) be the two faces whose boundaries contain \(e\). 
We denote by \(Q_e\) the quadrilateral \(vv_{f_1}wv_{f_2}\) in \(I(G)\).

\begin{figure}[htbp]
\centering
\includegraphics[scale=0.025]{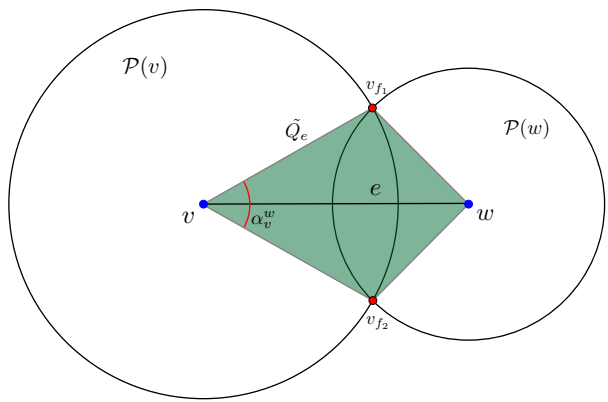}
\caption{Quadrilateral}
\label{5}
\end{figure}

Thus each edge of the original graph \(G\) corresponds uniquely to a quadrilateral in the incidence graph, as illustrated in Figure \ref{5}. 
From this data we construct a Euclidean quadrilateral \(\bar{Q}_e\) satisfying
\[
\angle vv_{f_1}w = \angle vv_{f_2}w = \pi-\Theta(e),\qquad 
|vv_{f_1}|=|vv_{f_2}|=r(v),\qquad 
|wv_{f_1}|=|wv_{f_2}|=r(w),
\]
as shown in Figure \ref{5}.

We write \(\alpha_{(e,v)}=\alpha_v^{w}\in(0,2\pi)\) for the angle \(\angle v_{f_1}vv_{f_2}\) in \(\tilde{Q}_e\) (see Figure \ref{5}), which is given explicitly by
$$
\alpha_v^{w}\bigl(r(v),r(w),\Theta(e)\bigr)=2\arccos\frac{r(v)+r(w)\cos\Theta(e)}{\sqrt{r(v)^2+r(w)^2+2r(v)r(w)\cos\Theta(e)}},\qquad 
(r(u),r(v))\in\mathbb{R}_+^2.
$$

Since this expression depends only on the ratio \(q_{vw}:=r(v)/r(w)\) and  \(\Theta(e)\), we may also write it as \(\alpha_v^{w}(q_{vw},\Theta(e))\). 
Gluing all the Euclidean quadrilaterals \(\bar{Q}_e\) along the edges of \(I(G)\) produces a piecewise flat metric \(g=g(r,\Theta)\) on \(\Sigma\) (for the gluing procedure we refer to~\cite[chapter~3]{BBS01}). 
The metric \(g\) is flat away from the points of \(V\cup V_F\); the vertices \(v\in V\) and the dual vertices \(v_f\in V_F\) may become cone points. 
In the present work we always assume that condition~\((C_1)\) holds for the pair \((\mathcal{D},\Theta)\), which guarantees that the points of \(V_F\) are not conical.

For a vertex \(v\in V\) we define its \textbf{cone angle} by
\[
\alpha_v=\sum_{e:\,v<e}\alpha_{(e,v)}=\sum_{w\sim v}\alpha_v^{w},
\]
and the \textbf{vertex curvature} (or \textbf{discrete Gauss curvature}) at \(v\) by
\[
K_v = 2\pi-\alpha_v.
\]

\begin{definition}
Let \(\mathcal{D}\) be a disk cellular decomposition and let \(\Theta\in(0,\pi)^E\) satisfy condition~\((C_1)\). 
A circle pattern metric \(r\) is said to support an \textbf{embedded planar ideal circle pattern} (embedded ICP) if the following two conditions hold:

\begin{itemize}
    \item[$(I)$] \(\alpha_v = 2\pi\) for every interior vertex \(v\in\operatorname{Int}(V)\).
    \item[$(II)$] There exists an isometric embedding \(\eta\colon(\Omega,g(\Theta,r))\to(\mathbb{R}^2,\mathrm{d}s^2)\), where \(\mathrm{d}s^2\) denotes the standard Euclidean metric.
\end{itemize}
\end{definition}

For brevity we shall not distinguish between the objects \(\mathcal{P}(v),\tilde{Q}_e\) and their images under \(\eta\). 
Observe that an ICP in \(\mathbb{R}^2\) as defined in Definition \ref{section-icp} is embedded precisely when the interiors of the quadrilaterals associated with different edges are mutually disjoint. In this paper, all ideal circle patterns are assumed to be \emph{embedded planar ideal circle patterns} unless otherwise specified.

In \cite{GeIIP1}, Ge-Yu-Zhou proved that for every IAG $(G,\Theta)$ with $\sup_{e\in E}\Theta(e)<\pi$, there exists an embedded ideal circle pattern ($\ICP$) that realizes it.
$\ICP$ embeddings can be classified into different types, such as parabolic and hyperbolic.
We say that the pair $(G,\Theta)$ is 
$\ICP$-parabolic if there exists an embedded ideal circle pattern ($\ICP$) realizing it that is locally finite in 
$\mathbb{R}^2$. 
Similarly, $(G,\Theta)$ is said to be $\ICP$-hyperbolic 
if there exists an $\ICP$ realizing it that is locally finite in the unit disk $\mathbb{D}$.

We recall the following rigidity results
from \cite{GeIIP1}. Related rigidity theorems
for infinite circle patterns were established
by B\"ucking \cite{Buck19}.

\begin{theorem}[Rigidity of parabolic ICPs]
Let $\mathcal{D}=(V, E, F)$ be an infinite disk cellular decomposition. Assume $\sup_{e\in E}\Theta(e)<\pi$ and conditions $(C_1)$ and $(C_2')$. If $\mathcal{P}$ and $\mathcal{P}^*$ are two embedded ICPs of $\mathcal{D}$, and $\mathcal{P}$ is locally finite in $\mathbb{R}^2$, then there exists an affine transformation $h$ such that $\mathcal{P}^*=h(\mathcal{P})$.
\end{theorem}

\begin{theorem}[Rigidity of hyperbolic ICPs]
Let $\mathcal{D}=(V, E, F)$ be an infinite disk cellular decomposition. Assume the intersection angle $\Theta \in(0, \pi)^E$ satisfies $(C_1)$ and $(C_2)$. If $\mathcal{P}$ and $\mathcal{P}^*$ are two ICPs of $\mathcal{D}$, and $\mathcal{P}$ is locally finite in the unit disk $\mathbb{D}$, then there exists a Mobius transformation $h$ such that $\mathcal{P}^*=h(\mathcal{P})$.
\end{theorem}

This implies that for an IAG $(G,\Theta)$ with $\sup_{e\in E}\Theta(e)<\pi$,
if it is $\ICP$-parabolic, then the $\ICP$ embedding is unique up to affine transformation;
if it is $\ICP$-hyperbolic, then the $\ICP$ embedding is unique up to Mobius transformation.

A basic feature of an embedded \(\ICP\) is that the radii of neighboring circles are comparable.
The local Ring Lemma gives such a comparison with a constant that may depend on the vertex and the neighborhood, while the uniform Ring Lemma yields bounds that depend only on global parameters. 
The following two lemmas state these estimates.

\begin{lemma}[Local Ring Lemma for $\ICP$ \cite{GeIIP1}]\label{lem:local-ring}
Let $\mathcal{P}$ be an embedded $\ICP$ of a finite disk cellular decomposition $\mathcal{D}=(V,E,F)$ with
intersection angles $\Theta$ satisfying
$\sup_{e\in E}\Theta(e)<\pi$ and conditions
{\rm ($C_1$)} and {\rm ($C_2$)}.
If $v\in V$ satisfies $B_{6\pi/\varepsilon_1}(v)\cap\partial V=\varnothing$,
then there exists a constant $C=C(v,\Theta,\mathcal{D})>0$ such that, for every edge $e=uv$,
\[ 
\frac{r(u)}{r(v)}\ \ge C(v,\Theta,\mathcal{D}).
\]
\end{lemma}

\begin{lemma}[Uniform Ring Lemma for $\ICP$ \cite{GeIIP1}]\label{lem:uniform-ring}
Let \(\mathcal P\) be an embedded \(\ICP\) of a finite disk cellular decomposition
\(\mathcal{D}=(V,E,F)\) with vertex degree bounded by \(N\)
and intersection angles \(\Theta\) satisfying
\(\sup_{e\in E}\Theta(e)<\pi\), \((C_1)\), and \((C_2')\).
If \(v\in V\) satisfies
\(
B_{6\pi/\varepsilon_1}(v)\cap \partial V=\varnothing,
\)
then there exist a constant $C=C(N,\varepsilon_1,\varepsilon_0)$
such that for every edge $e=uv$,
\[
\frac{r(u)}{r(v)}\ \ge\ C(N,\varepsilon_1,\varepsilon_0)
\]
\end{lemma}

In this paper, we prove the following refined ring lemma for $\ICP$ in Section~\ref{sec:refined}.

\begin{theorem}[Refined Ring Lemma for $\ICP$]
    Let $(G, \Theta)$ be a tame $\IIA$ with $\Theta \in [0, \pi - \varepsilon_1]^E$. Let its embedding on the plane have radius $r: V \rightarrow \mathbb{R}$. Then there exists a uniform constant $C = C(\varepsilon_0, \varepsilon_1)$ such that, for any edge $u \sim v$:

$$\frac{r(v)}{r(u)} > e^{-C\cdot S(u)}$$

where $S(u)=\sum_{v' \sim u} \deg(v')$, called the flower degree of $u$.
\end{theorem}

\subsection{Infinite ideal polyhedra in $\mathbb{H}^3$}
\label{section-iip-H3}
\begin{definition}[Infinite Ideal Polyhedron]
	An infinite ideal polyhedron $(\IIP)$ in hyperbolic three-space $\mathbb{H}^3$ 
	is defined as the convex hull of a countable collection of isolated points located on the sphere at infinity, 
	which are referred to as ideal vertices.
\end{definition}

Let $P$ be such an infinite ideal polyhedron in $\mathbb{H}^3$. 
	We denote by $\partial_0 P$ the boundary of $P$ inside $\mathbb{H}^3$, 
	and by $\partial_\infty P$ the ideal boundary of $P$ in the Poincar\'e ball model. 
	It is clear that $\partial_\infty P$ coincides with the set of isolated ideal vertices forming $P$.

	In \cite{GeIIP1}, Ge-Yu-Zhou proved that when $\sup_{e\in E}\Theta(e)<\pi$ and $(C_1)$,$(C_2)$ hold, there exists an infinite ideal polyhedron (abbreviated $\IIP$) in $\mathbb{H}^3$ that is combinatorially equivalent to the Poincar\'e dual of $\mathcal{D}$ with the dihedral angle.
	$\Theta_P(e)=\Theta(e^{*})$. This implies that for every $\IIA$ with $\sup_{e\in E}\Theta(e)<\pi$, there exists an $\IIP$ in $\mathbb{H}^3$ that is combinatorially equivalent to the Poincar\'e dual of $\mathcal{D}$ with the dihedral angle.
	
	Since $\mathbb{S}^2$ is compact, the set $\partial_\infty P$ must possess accumulation points on $\mathbb{S}^2$. 
	Consequently, the union $\partial_0 P \cup \partial_\infty P$ is not homeomorphic to the sphere, 
	but rather to the sphere with certain limit point sets removed. 
	The classification of $P$ into parabolic or hyperbolic type 
	reflects the geometric structure of these accumulation points.
    
\begin{definition}[Parabolic and Hyperbolic Types]\label{def:comfor}
	Let $P$ be an infinite ideal polyhedron in $\mathbb{H}^3$.
	\begin{enumerate}[(1)]
		\item $P$ is said to be of parabolic type (abbreviated PIIP) 
		if $\partial_\infty P$ has exactly one accumulation point on $\mathbb{S}^2$.
		\item $P$ is said to be of hyperbolic type (abbreviated HIIP) 
		if the set of accumulation points of $\partial_\infty P$ forms a circle 
		$C_\infty \subset \mathbb{S}^2$, and all points of $\partial_\infty P$ 
		lie on the same side of $C_\infty$.
	\end{enumerate}
\end{definition}

Ge-Yu-Zhou proved that for an IAG $(G,\Theta)$ with $\Theta\in (0,\pi)^{E}$,
all of the corresponding PIIPs are isometric;
for a tame IAG $(G,\Theta)$ with $\sup_{e\in E}\Theta(e)<\pi$,
all of the corresponding HIIPs are isometric;
(\cite[Thm.~1.7]{GeIIP1}).

\subsection{Unimodular random rooted graphs}

\medskip

A rooted graph is a pair $(G,\rho)$ consisting of a connected, locally finite graph $G$
and a distinguished vertex $\rho \in V(G)$, called the root.
Two rooted graphs $(G,\rho)$ and $(G',\rho')$ are said to be isomorphic
if there exists a graph isomorphism $\varphi:G\to G'$ satisfying $\varphi(\rho)=\rho'$.
Let $\mathcal{G}_\bullet$ denote the space of all isomorphism classes of rooted, locally finite graphs,
equipped with the Benjamini-Schramm topology (local weak topology).

A random rooted graph is a random variable taking values in $\mathcal{G}_\bullet$.
Such a graph $(G,\rho)$ is said to be unimodular
if it satisfies the Mass-Transport Principle (MTP):
for every nonnegative Borel measurable function 
$f:\mathcal{G}_{\bullet\bullet}\to[0,\infty]$ 
depending on a graph with an ordered pair of distinguished vertices $(\rho,v)$,
\[
\mathbb{E}\!\left[\sum_{v\in V(G)} f(G,\rho,v)\right]
\;=\;
\mathbb{E}\!\left[\sum_{v\in V(G)} f(G,v,\rho)\right],
\]
where $\mathcal{G}_{\bullet\bullet}$ denotes the space of doubly rooted graphs.

A map is a proper embedding of a connected graph into a surface considered up to orientation-preserving homeomorphism. 
In particular, a planar map is an embedding of a graph into the sphere $\mathbb{S}^2$ 
(or equivalently into the plane $\mathbb{R}^2$)
When a map is endowed with a distinguished vertex $\rho$, 
we refer to the pair $(M,\rho)$ as a rooted map.

A random rooted map $(M,\rho)$ is said to be
unimodular if the Mass-Transport Principle
holds for all nonnegative Borel functions
on doubly rooted maps.
For a random rooted IAG $(G,\rho,\Theta)$,
unimodularity is understood in the sense
of angle-marked graphs: all isomorphisms
and mass transports preserve the angle marks.

\subsection{Vertex extremal length and random walks}

The vertex extremal length (VEL), introduced by Schramm~\cite{Schramm1993}, 
is a discrete analogue of the classical extremal length on Riemann surfaces, 
and provides a powerful criterion for distinguishing recurrence and transience 
of random walks on planar graphs.

Let $G=(V,E)$ be an infinite connected graph, and let $A\subset V$ be a finite subset. 
Denote by $\Gamma(A,\infty)$ the family of all infinite paths in $G$ starting from $A$ 
that cannot be contained in any finite subset of $V$. 
A function $m:V\to[0,\infty)$ is said to be $\Gamma(A,\infty)$-admissible if
\[
\sum_{v\in\gamma} m(v) \ge 1, \qquad \text{for every } \gamma\in\Gamma(A,\infty).
\]
The vertex extremal length of $\Gamma(A,\infty)$ is defined as
\[
\mathrm{VEL}(A,\infty)
= \sup_{m} \frac{1}{\sum_{v\in V} m(v)^2},
\]
where the supremum is taken over all $\Gamma(A,\infty)$-admissible functions $m$.

\begin{definition}[VEL-parabolic and VEL-hyperbolic graphs]
	An infinite graph $G=(V,E)$ is called VEL-parabolic if there exists a finite subset 
	$A\subset V$ such that $\mathrm{VEL}(A,\infty)=\infty$;  
	otherwise, $G$ is called VEL-hyperbolic.
\end{definition}

\medskip

Let $G$ be a locally finite (possibly weighted) graph.
The simple random walk $(X_n)_{n\ge0}$ on $G$ is the Markov chain with transition
probabilities
\[
p(x,y) = \frac{w(x,y)}{w(x)}, \qquad
w(x) = \sum_{e\ni x} w(e),
\]
where $w:E(G)\to\mathbb{R}_+$ is the edge-weight function (for unweighted graphs, $w\equiv1$).
We allow infinite degree provided $w(x)<\infty$ for all $x\in V(G)$.

\medskip

Denote by $\mathcal{G}^{\leftrightarrow}$ (resp.\ $\mathcal{M}^{\leftrightarrow}$)
the space of isomorphism classes of graphs (resp.\ maps) endowed with a
bi-infinite path $(x_n)_{n\in\mathbb{Z}}$,
equipped with the local topology.  
Given a random rooted graph $(G,\rho)$, let $(X_n)_{n\ge0}$ and $(X_{-n})_{n\ge0}$ 
be two independent simple random walks started from $\rho$, and
consider $(G,(X_n)_{n\in\mathbb{Z}})$ as a random element of $\mathcal{G}^{\leftrightarrow}$.

\begin{definition}[Stationarity and reversibility]
	A random rooted graph $(G,\rho)$ is stationary if
	\[
	(G,\rho)\ \stackrel{d}{=}\ (G,X_1),
	\]
	and reversible if
	\[
	(G,\rho,X_1)\ \stackrel{d}{=}\ (G,X_1,\rho),
	\]
	as doubly rooted graphs.  
	Equivalently, $(G,\rho)$ is reversible if and only if
	\[
	(G,(X_n)_{n\in\mathbb{Z}})\ \stackrel{d}{=}\ (G,(X_{n+k})_{n\in\mathbb{Z}}), \qquad \forall\,k\in\mathbb{Z},
	\]
	that is, the bi-infinite random walk is stationary under the shift.
\end{definition}

\medskip

Reversibility and unimodularity are closely related.
If $(G,\rho)$ is reversible, then reweighting its law by $\deg(\rho)^{-1}$, i.e.\
by the Radon-Nikodym derivative
\[
\frac{\deg(\rho)^{-1}}{\mathbb{E}[\deg(\rho)^{-1}]},
\]
yields an equivalent unimodular random rooted graph.
Conversely, if $(G,\rho)$ is unimodular with $\mathbb{E}[\deg(\rho)]<\infty$,
biasing by $\deg(\rho)$ produces a reversible random rooted graph.
Thus, reversible and unimodular laws are in one-to-one correspondence
when the expected degree is finite.

Every VEL-hyperbolic graph is transient for simple random walk. Moreover, for graphs of bounded degree the converse holds: transience is equivalent to VEL hyperbolicity \cite{HeSc95}.

\medskip

An event $A\subset\mathcal{G}^{\leftrightarrow}$ is shift-invariant
if $(G,(X_n)_{n\in\mathbb{Z}})\in A$ implies 
$(G,(X_{n+k})_{n\in\mathbb{Z}})\in A$ for all $k\in\mathbb{Z}$.
A reversible or unimodular random graph is said to be ergodic
if every invariant event has probability either $0$ or $1$.

\begin{theorem}[Characterisation of ergodicity {\cite{AldLy07} or \cite[Theorem 3.1]{AHNR16}}]
	Let $(G,\rho)$ be a unimodular random rooted graph with $\mathbb{E}[\deg(\rho)]<\infty$
	(or a reversible random rooted graph).  
	The following are equivalent:
	\begin{enumerate}[(1)]
		\item $(G,\rho)$ is ergodic;
		\item every rerooting-invariant event $A\subset\mathcal{G}_\bullet$
		satisfies $\mathbb{P}(A)\in\{0,1\}$;
		\item the law of $(G,\rho)$ is an extreme point of the weakly closed convex set
		of unimodular (resp.\ reversible) laws.
	\end{enumerate}
\end{theorem}

\medskip

By Choquet's theorem, every unimodular random rooted graph admits an
ergodic decomposition: it can be represented as a mixture of ergodic
unimodular random rooted graphs.  
Hence, to establish almost-sure properties for a general unimodular law,
it suffices to consider the ergodic case.

\subsection{Invariant amenability}

We briefly recall the notion of invariant amenability, following
Aldous and Lyons~\cite{AldLy07} and
Angel et al. \cite{AHNR16,map}.

\medskip
A weighted graph is a graph equipped with a weight function \(w: E \to \mathbb{R}_+\). This function is extended to vertices by \(w(x) = \sum_{e \ni x} w(e)\), and---allowing a mild notational abuse---to subsets of edges or vertices by additivity. For an infinite weighted graph, the edge Cheeger constant is defined as  

\[
\mathbf{i}_E(G)=\inf\Bigl\{\frac{w(\partial_E W)}{w(W)} : \emptyset \neq W \subset V,\ W \text{ finite}\Bigr\},
\]  
where \(\partial_E W\) denotes the set of edges having precisely one endpoint in \(W\). The graph is called \textit{amenable} if \(\mathbf{i}_E(G)=0\), and \textit{non-amenable} if \(\mathbf{i}_E(G)>0\).

Let $(G,\rho)$ be a unimodular random rooted graph.
A percolation on $(G,\rho)$ is a random labeling 
$\omega:E(G)\cup V(G)\to\{0,1\}$ such that the marked graph $(G,\rho,\omega)$
is unimodular.  We think of $\omega$ as specifying the open edges and vertices, 
and assume that whenever an edge is open, both its endpoints are open as well.
The cluster of a vertex $v$ in $\omega$ is
\[
K_\omega(v) := \{u\in V(G) : u \text{ is connected to } v
\text{ by a path of open edges in } \omega \}.
\]
A percolation is called finitary if every cluster $K_\omega(v)$ is finite almost surely.

\medskip

The invariant Cheeger constant of an ergodic unimodular random rooted graph 
$(G,\rho)$ is defined by
\begin{equation}\label{eq:iinv}
	i_{\mathrm{inv}}(G,\rho)
	:= \inf_{\omega}
	\mathbb{E}\!\left[
	\frac{|\partial_E K_\omega(\rho)|}{|K_\omega(\rho)|}
	\right],
\end{equation}
where the infimum is taken over all finitary percolations $\omega$ on $(G,\rho)$,
and $\partial_E K_\omega(\rho)$ denotes the set of edges having exactly one endpoint in $K_\omega(\rho)$.

\medskip

Associated to $(G,\rho)$ is another quantity describing the mean internal degree
in finitary percolations:
\[
\alpha(G,\rho)
:= \sup\bigl\{
\mathbb{E}[\deg_\omega(\rho)] :
\omega \text{ a finitary percolation on } (G,\rho)
\bigr\},
\]
where $\deg_\omega(\rho)$ is the degree of $\rho$ in the subgraph induced by $\omega$
(and $\deg_\omega(\rho)=0$ if $\rho\notin\omega$).
By the mass-transport principle (see~\cite[Proposition 8.3]{AldLy07}), one has
\[
\mathbb{E}[\deg(\rho)]
= i_{\mathrm{inv}}(G,\rho) + \alpha(G,\rho).
\]
In particular, when $\mathbb{E}[\deg(\rho)]<\infty$,
the invariant Cheeger constant $\iota_{\mathrm{inv}}(G,\rho)$ is positive
if and only if $\alpha(G,\rho) < \mathbb{E}[\deg(\rho)]$.

\begin{definition}[Invariant amenability]
	An ergodic unimodular random rooted graph $(G,\rho)$ is called
	invariantly amenable if $i_{\mathrm{inv}}(G,\rho)=0$,
	and invariantly non-amenable otherwise.
\end{definition}

\subsection{Embeddings and boundaries}

We now review several geometric embeddings of planar graphs that will be
relevant to our discussion\textemdash namely, the circle packing (CP) embedding,
the ideal circle pattern ($\ICP$) embedding.  Each of these constructions provides a concrete model for the boundary at infinity of the corresponding random walk or harmonic functions.

\medskip
\paragraph{\textbf{Circle packing (CP) embedding.}}
Let $G=(V,E)$ be a locally finite planar triangulation.
By the Koebe-Andreev-Thurston circle packing theorem,
there exists a circle packing $\mathcal{P}=\{\mathcal{P}(v)\}_{v\in V}$
in the Riemann sphere $\widehat{\mathbb{C}}$ whose tangency graph is $G$.
The embedding is unique up to Möbius transformations.
Depending on the conformal type of $\mathcal{P}$, we distinguish two cases:

In the hyperbolic case, the circles accumulate at the unit circle $\partial\mathbb{D}$,
which serves as a geometric compactification of $G$.
\cite{ABGN16} proved that $\partial\mathbb{D}$
simultaneously realizes the Poisson and Martin boundaries
of the simple random walk on the triangulation $G$ with bounded degree.

\medskip
\paragraph{\textbf{Ideal circle pattern ($\ICP$) embedding.}}
The circle packing model can be generalized to allow prescribed intersection angles.
Given a cellular decomposition $\mathcal{D}=(V,E,F)$ and a function
$\Theta:E\to(0,\pi)$ satisfying the combinatorial angle conditions $(C_1)$ and $(C_2')$,
there exists an ideal circle pattern 
$\mathcal{P}=\{\mathcal{P}(v)\}_{v\in V}$
realizing $(\mathcal{D},\Theta)$.
Each face $f$ corresponds to an intersection point $v_f$.

In the hyperbolic case, the intersection points accumulate at $\partial\mathbb{D}$,
which again coincides with the geometric and the Poisson boundaries.
The $\ICP$ embedding extends Thurston’s discrete conformal theory to
configurations with non-tangential intersections,
and corresponds bijectively to ideal hyperbolic polyhedra
whose dihedral angles are $\Theta(e)$
(see~\cite{GeIIP2,GeIIP1}).

These constructions relate geometric embeddings to harmonic functions
under their respective hypotheses. In the present paper,
Theorem~\ref{thm:icp-poisson} identifies the unit circle, equipped with the
exit measures, as a realization of the Poisson boundary for the
unimodular random ideal angled graphs considered there.

\section{Dichotomy theorem}

\subsection{Proof of Theorem 1.1}

\begin{theorem}
\label{thm:section3-1}
Let $(G, \rho, \Theta)$ be a unimodular random rooted $\IIA$, then
\(
\mathbb{E}[T(\rho)]=\mathbb{E}[\kappa(\rho)].
\)
If further assume that $\mathbb{E}[\deg(\rho)]<\infty$, then
\(
\mathbb{E}[\kappa(\rho)], \mathbb{E}[T(\rho)]
\)
are finite.
\end{theorem}

\begin{proof}
For a face $f$, 
\[
\partial f = (v_1,e_1,v_2,e_2,\dots,v_{d},e_{d},v_1),
\qquad d=\deg(f),
\]

For a vertex $v$ and a face $f$, define the boundary multiplicity
\begin{equation}\label{eq:kfv}
k_f(v):=\#\{\, i\in\{1,\dots,d\} \,:\, v_i=v \,\},
\end{equation}
namely the number of times $v$ appears on $\partial f$. For an $\IIA$, since it is a disk cellular decomposition, each face boundary is a simple cycle; hence
\(
k_f(v)\equiv 1 \ \text{or} \ 0 
\)
.

Since the underlying graph is simple, every edge $e$ has two distinct endpoints, say
$\partial e=\{x,y\}$ with $x\neq y$, and hence $\sum_{u\in\partial e}1=2$.

\begin{figure}[htbp]
\centering
\includegraphics[scale=0.018]{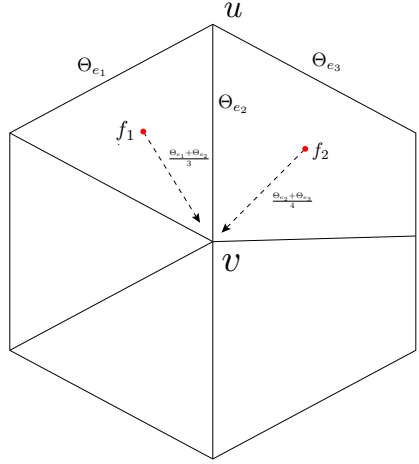}
\caption{Mass transport $m(u,v)$}
\label{6}
\end{figure}

Define a mass transport $m=m(u,v)$ by (Figure~\ref{6})
\begin{equation}\label{eq:mass-transport-def}
m(u,v)
:=\sum_{f \in F(G)}\;\sum_{i=1}^{\deg(f)} \mathbf{1}_{\{v_i=v\}}
\;\sum_{e\in\partial f}\;\mathbf{1}_{\{u\in\partial e\}}
\frac{\Theta_e}{2\,\deg(f)}.
\end{equation}
Equivalently: for each face $f$, for each boundary edge $e\in\partial f$, each endpoint
of $e$ sends mass $\Theta_e/(2\deg(f))$ to each boundary-vertex occurrence of $f$; then
$m(u,v)$ aggregates over all faces and all occurrences where that boundary vertex equals $v$.
In particular, $m(u,v)\ge 0$. 
The definition is measurable and isomorphism-invariant since it depends only on the rooted
graph structure and angles $\{\Theta_e\}$.

Now, we compute the total mass sent from a fixed vertex.
Fix $u\in V(G)$. We compute $\sum_{v\in V(G)} m(u,v)$.

Starting from \eqref{eq:mass-transport-def}, since $G$ is locally finite and each face has finite degree, for fixed $u$ only finitely many
terms in the sums are nonzero; hence all sums below are finite and we may freely change the order.
\begin{align*}
\sum_{v\in V(G)} m(u,v)
&=
\sum_{v\in V(G)}
\sum_{f\in F(G)}\sum_{i=1}^{\deg(f)} \mathbf{1}_{\{v_i=v\}}
\sum_{e\in\partial f}\mathbf{1}_{\{u\in\partial e\}}\frac{\Theta_e}{2\,\deg(f)}
\\
&=
\sum_{f\in F(G)}\sum_{i=1}^{\deg(f)}
\sum_{e\in\partial f}\mathbf{1}_{\{u\in\partial e\}}\frac{\Theta_e}{2\,\deg(f)}
\sum_{v\in V(G)}\mathbf{1}_{\{v_i=v\}}
\\
&=
\sum_{f\in F(G)}\sum_{i=1}^{\deg(f)}
\sum_{e\in\partial f}\mathbf{1}_{\{u\in\partial e\}}\frac{\Theta_e}{2\,\deg(f)}
\\
&=
\sum_{f\in F(G)} \deg(f)\cdot
\sum_{e\in\partial f}\mathbf{1}_{\{u\in\partial e\}}\frac{\Theta_e}{2\,\deg(f)}
\\
&=
\sum_{f\in F(G)}\sum_{e\in\partial f}\mathbf{1}_{\{u\in\partial e\}}\frac{\Theta_e}{2}.
\end{align*}

Fix an edge $e$ incident to $u$. Since the underlying graph is simple, exactly one endpoint
of $e$ equals $u$, hence $\mathbf{1}_{\{u\in\partial e\}}=1$. Moreover, $e$ has exactly two
face-incidences counted with multiplicity, so it appears once in the boundary walk of each
incident face-incidence. Therefore $e$ contributes $\Theta_e/2$ for each face-incidence and
hence contributes $\Theta_e$ in total. Summing over all edges $e \ni u$ yields
\[
\sum_{v\in V(G)} m(u,v)=\sum_{e \ni u}\Theta_e =2\pi-T(u).
\]
In particular,
\begin{equation}\label{eq:outgoing-equals-T}
\sum_{v\in V(G)} m(\rho,v)=2\pi-T(\rho)\qquad\text{a.s.}
\end{equation}

Now, we compute the total mass received by a fixed vertex. Fix $v\in V(G)$ and compute $\sum_{u\in V(G)} m(u,v)$.

From \eqref{eq:mass-transport-def},
\begin{align*}
\sum_{u\in V(G)} m(u,v)
&=
\sum_{u\in V(G)}\sum_{f\in F(G)}\sum_{i=1}^{\deg(f)} \mathbf{1}_{\{v_i=v\}}
\sum_{e\in\partial f}\mathbf{1}_{\{u\in\partial e\}}\frac{\Theta_e}{2\,\deg(f)}
\\
&=
\sum_{f\in F(G)}\sum_{i=1}^{\deg(f)} \mathbf{1}_{\{v_i=v\}}
\sum_{e\in\partial f}\left(\sum_{u\in V(G)}\mathbf{1}_{\{u\in\partial e\}}\right)\frac{\Theta_e}{2\,\deg(f)}.
\end{align*}
For each edge $e$ in a simple graph, $\partial e$ has exactly two endpoints, hence
\[
\sum_{u\in V(G)}\mathbf{1}_{\{u\in\partial e\}}=2.
\]
Therefore
\begin{align*}
\sum_{u\in V(G)} m(u,v)
&=
\sum_{f\in F(G)}\sum_{i=1}^{\deg(f)} \mathbf{1}_{\{v_i=v\}}
\sum_{e\in\partial f}\frac{\Theta_e}{\deg(f)}
\\
&=
\sum_{f\in F(G)}\sum_{i=1}^{\deg(f)} \mathbf{1}_{\{v_i=v\}}
\frac{1}{\deg(f)}\sum_{e\in\partial f}\Theta_e.
\end{align*}
By ($C_1$) condition in the definition of $\IIA$, $\sum_{e\in\partial f}\Theta_e=(\deg(f)-2)\pi$, hence
\begin{align*}
\sum_{u\in V(G)} m(u,v)
&=
\sum_{f\in F(G)}\sum_{i=1}^{\deg(f)} \mathbf{1}_{\{v_i=v\}}
\frac{(\deg(f)-2)\pi}{\deg(f)}
\\
&=
\sum_{f\in F(G)} k_f(v)\,\frac{(\deg(f)-2)\pi}{\deg(f)}.
\\
&=
\sum_{f\in F(G)} \frac{(\deg(f)-2)\pi}{\deg(f)}.
\end{align*}

The last expression equals $\theta(v)$. In particular,
\begin{equation}\label{eq:incoming-equals-theta}
\sum_{u\in V(G)} m(u,\rho)=\theta(\rho)\qquad\text{a.s.}
\end{equation}

Now, we check the integrability.
Since $\Theta_e\in(0,\pi)$,
\[
T(\rho)=2\pi-\sum_{e\sim \rho}\Theta_e \ge 2\pi-\pi\,\deg(\rho),
\]
so $\mathbb{E}[T(\rho)]>-\infty$ follows from $\mathbb{E}[\deg(\rho)]<\infty$.

Also $0\le \frac{(\deg(f)-2)\pi}{\deg(f)}\le \pi$ for every face $f$, and the total number of
face corners incident to $\rho$ equals $\deg(\rho)$ in a simple graph. Hence
\[
\theta(\rho)=\sum_{f\ni \rho}\frac{(\deg(f)-2)\pi}{\deg(f)}\le \pi\,\deg(\rho),
\]
so $\mathbb{E}[\kappa(\rho)]<\infty$.

Finally, we apply the Mass Transport Principle.
Because $(G,\rho,\Theta)$ is unimodular and $m\ge 0$ is measurable and isomorphism-invariant,
the Mass Transport Principle gives
\[
\mathbb{E}\Big[\sum_{v\in V(G)} m(\rho,v)\Big]
=
\mathbb{E}\Big[\sum_{u\in V(G)} m(u,\rho)\Big].
\]
Using \eqref{eq:outgoing-equals-T} and \eqref{eq:incoming-equals-theta}, the left-hand side
equals $2\pi-\mathbb{E}[T(\rho)]$ and the right-hand side equals $\mathbb{E}[\theta(\rho)]$.
Therefore,
\[
\mathbb{E}[T(\rho)]=\mathbb{E}[\kappa(\rho)].
\]
\end{proof}

\subsection{Proof of Theorem \ref{thm:ihp-dichotomy} and \ref{thm:iag-dichotomy}.}
By Theorem \ref{thm:ihp-gb}, the unimodular Gauss-Bonnet identity
takes the equivalent graph-theoretic form $\mathbb E[T(\rho)]=\mathbb E[\kappa(\rho)]$. The curvature dichotomy for unimodular random planar
maps established by Angel-Hutchcroft-Nachmias-Ray
in \cite{map} implies that
\(\mathbb E[\kappa(\rho)]\leq 0\), with equality if
and only if \(G\) is invariantly amenable, or
equivalently, almost surely VEL-parabolic.
Consequently, statements (1)-(3) of
Theorem \ref{thm:iag-dichotomy} are equivalent.

Suppose $(G, \rho,\Theta)$ is $\ICP$-hyperbolic, and let $\mathcal{P}$ be an ideal circle pattern of $G$. Then, we embed $G$ into $\mathbb{D}$ and connect the hyperbolic centers of the circles in $\mathcal{P}$ with hyperbolic geodesics, thereby representing each face of $G$ that is a triangle as a hyperbolic triangle. By the rigidity theorem(see \cite{GeIIP1}), this embedding is unique up to hyperbolic isometries, depending only on the isomorphism type of the IAG
$(G,\Theta)$.

Similarly, if $(G, \rho,\Theta)$ is $\ICP$-parabolic with $\sup_e\Theta(e)<\pi$, we embed it into $\mathbb{C}$ by connecting the Euclidean centers of the circles with straight lines.

Then, for each face of $(G,\rho)$, we connect $v_f$ to all vertices of that face by geodesics and get a new graph $\tilde{G}$. Then, we define a mass transport as follows. For each edge $(u,v)$ in $G$, there are two faces $f_1=(u,v,...)$ and $f_2=(u,v,...)$ containing this edge and transport $\angle v_{f_1}uv_{f_2}$ from $u$ to $v$. Besides, the transport from $u$ to itself has a term for each face containing $u$. Since these angles are independent of the choice of circle pattern by rigidity(see \cite{GeIIP1}), we know that this mass transport is well-defined.

\begin{figure}[htbp]
\centering
\includegraphics[scale=0.026]{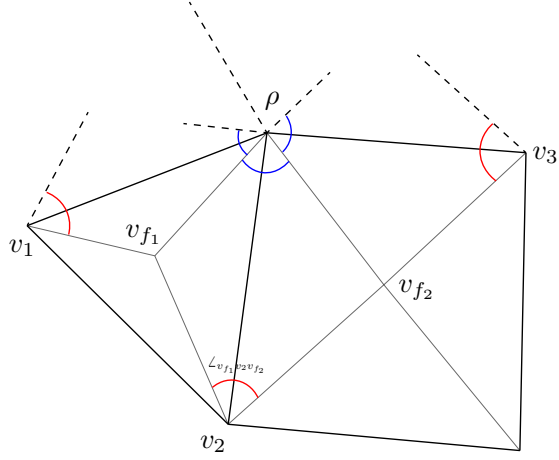}
\caption{Mass transport between $\rho$ and $v_2$}
\label{7}
\end{figure}

Specifically, we define the mass transport $F: \mathcal{G}_{\bullet \bullet} \rightarrow \mathbb{R}_{+}$ by (Figure~\ref{7})
\begin{equation}
	F(G, u, v)= \begin{cases}\angle v_{f_1}uv_{f_2} & \text { if } v \sim u, v \neq u, \\ 
		2\pi & \text { if } v=u,\\
		0 & \text { otherwise.}\end{cases}
\end{equation}
Suppose all faces containing the root $\rho$ are $f_i=(v_i,\rho,v_{i+1},...)$, where $i=1, \ldots, n$. In particular, we define $v_{n+1}$ to be $v_1$ and $f_{n+1}$ to be $f_1$. Then, we have 
\begin{equation}
	\sum_{v} F(G, \rho, v)=4\pi,
\end{equation}
and
\begin{equation}
	\begin{aligned}
		\sum_{v} F(G, v, \rho) & =F(G, \rho, \rho)+\sum_{\substack{v \sim \rho, v \neq \rho}} F(G, v, \rho) \\
		& =2 \pi+\sum_{i=1}^n \angle v_{f_i} v_i v_{f_{i+1}} \\
		& =\sum_{i=1}^n \angle v_{f_i} \rho v_{f_{i+1}}+\sum_{i=1}^n\angle v_{f_i} v_i v_{f_{i+1}} \\
		& =2 \sum_{i=1}^n\left(\theta_i-\pi+\Theta_{\rho v_i}\right) \\
		& =4\pi-2T(\rho)+2 \sum_{i=1}^n\left(\theta_i-\pi\right).
	\end{aligned}
\end{equation}
Here \(\theta_i\) denotes the sum of the internal angles of the triangle \(\triangle(\rho,v_i,v_{f_i})\) in the drawing. In the Euclidean case, \(\theta_i=\pi\), whereas in the hyperbolic case, \(\theta_i<\pi\). Hence, we know that 
\begin{equation}
	\sum_{v} F(G, v, \rho)=4\pi-2T(\rho),
\end{equation}
when $(G,\rho,\Theta)$ is $\ICP$-parabolic and
\begin{equation}
	\sum_{v} F(G, v, \rho)<4\pi-2T(\rho),
\end{equation}
when $(G,\rho,\Theta)$ is $\ICP$-hyperbolic. Then, by the mass transport principle, we have $\mathbb{E}[T(\rho)]=0$ when $(G,\rho,\Theta)$ is $\ICP$-parabolic and $\mathbb{E}[T(\rho)]<0$ when $(G,\rho,\Theta)$ is $\ICP$-hyperbolic. 

Thus, if $(G,\rho,\Theta)$ is $\ICP$-hyperbolic, by The Dichotomy Theorem in \cite{map}, we know that $G$ is invariantly nonamenable and almost surely VEL-hyperbolic. Since a simple random walk on any VEL-hyperbolic graph is transient, it follows that $G$ is almost surely transient.

If $\sup_e\Theta(e)<\pi$, we know the ICP type of $(G,\rho,\Theta)$ has rigidity \cite{Buck19,GeIIP1}. Since $(G,\rho,\Theta)$ is ergodic and the $\ICP$ type does not depend on the choice of root, we know that the $\ICP$ type of $(G,\rho,\Theta)$ is almost surely constant. Then, from the previous proof, we know $E[T(\rho)] \leq 0$, and $E[T(\rho)] = 0$ if and only if $(G,\rho,\Theta)$ is $\ICP$-parabolic. Thus, by The Dichotomy Theorem in \cite{map}, we know that statements (1), (2), (3) and (4) are equivalent. Finally, the statements (3) and (5) are equivalent under the bounded degree assumption (see \cite{HeSc95}).

\begin{remark}[An alternative argument for invariant amenability]
\label{rem:T-amenable1}
Assume $\mathbb{E}[T(\rho)]=0$. By Theorem \ref{thm:ihp-gb} we have
$\mathbb{E}[T(\rho)]=\mathbb{E}[\kappa(\rho)]$, so
\[
\mathbb{K}(G,\rho):=\mathbb{E}\bigl[\kappa(\rho)\bigl]
\]
vanishes. 

Let $F$ be a sample from the free uniform spanning forest (FUSF) of $G$.
Angel-Hutchcroft-Nachmias-Ray \cite[Theorem~5.11 and Corollary~5.12]{map} (see also~\cite[\S5.5]{map})
give the identity
\[
\mathbb{E}[\deg_F(\rho)]
=\frac{1}{\pi}\,\mathbb{E}[\theta(\rho)]
=2-\frac{1}{\pi}\,\mathbb{K}(G,\rho).
\]
Since $\mathbb{K}(G,\rho)=0$, we get $\mathbb{E}[\deg_F(\rho)]=2$.
On the other hand, for the wired USF one always has
$\mathbb{E}[\deg_{WUSF}(\rho)]=2$ on unimodular random graphs~\cite[Proposition 7.1]{AldLy07},
and FUSF stochastically dominates WUSF; therefore equality of expectations forces
$FUSF=WUSF$. By the Aldous-Lyons amenability criterion (as used in~\cite[\S5.5]{map}),
this is equivalent to $(G,\rho)$ being invariantly amenable.
\end{remark}

Finally, Theorem \ref{thm:ihp-dichotomy} follows from the
polyhedral interpretation of ICP type
under the IHP/IAG duality. 

\subsection{The degree-six dichotomy for unimodular
triangulations}

We now consider the special case of
Theorem \ref{thm:iag-dichotomy} in which
\(\sup_{e\in E(G)}\Theta(e)\leq \pi/3\).
By the ideal-vertex condition (C1), every face of
the underlying cellular decomposition must be a
triangle, and \(\Theta(e)=\pi/3\) for every edge.
Consequently,
\[
T(\rho)
=2\pi-\sum_{v\sim\rho}\Theta_{\rho v}
=\frac{\pi}{3}(6-\deg(\rho)),
\]
and hence
\[
\mathbb E[T(\rho)]
=\frac{\pi}{3}
 \bigl(6-\mathbb E[\deg(\rho)]\bigr).
\]

Since \(\mathbb E[T(\rho)]\leq0\) by
Theorem \ref{thm:iag-dichotomy}, we obtain
\[
\mathbb E[\deg(\rho)]\geq6.
\]
Moreover, Theorem \ref{thm:iag-dichotomy} implies that
\[
\mathbb E[\deg(\rho)]=6
\quad\Longleftrightarrow\quad
G \text{ is almost surely ICP-parabolic},
\]
whereas
\[
\mathbb E[\deg(\rho)]>6
\quad\Longleftrightarrow\quad
G \text{ is almost surely ICP-hyperbolic}.
\]
The corresponding characterizations in terms
of invariant amenability and vertex extremal
length follow from the same theorem.

Thus, within the class of triangulations arising
from our ideal-angled setting, we recover the
degree-six criterion of
Angel-Hutchcroft-Nachmias-Ray \cite[Theorem 1.1]{AHNR16},
now expressed in terms of ideal circle-pattern
geometry.

\section{Boundary theory}

\subsection{Refined ring lemma for $\ICP$}\label{sec:refined}
Before studying the properties of random walks on $\IIA$, we first need to establish more precise geometric properties of the associated $\ICP$.
In this section we prove a refined ring lemma for $\ICP$ of tame $\IIA$.
Compared with the local ring lemma, we quantify the radius ratio by degrees:
for a vertex $v$ with neighbors $v_i$ and $d_i=\deg(v_i)$, the ratio is
controlled by an exponential factor in $\sum_i d_i$.
Compared with the uniform ring lemma, we weaken the conditions by removing the bounded degree assumption.

\begin{figure}[htbp]
\centering
\includegraphics[scale=0.025]{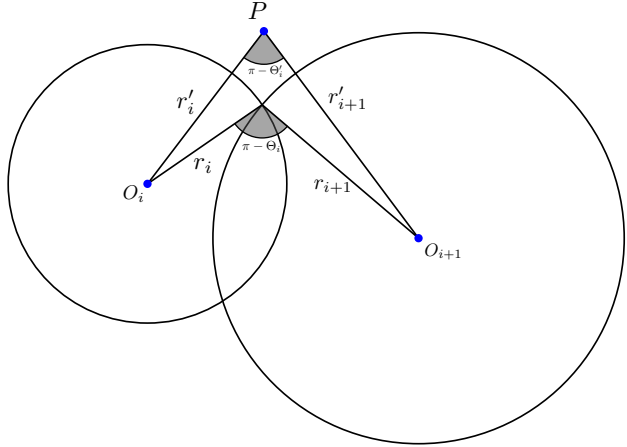}
\caption{Graph for Lemma~\ref{lemma1}}
\label{8}
\end{figure}

\begin{lemma}\label{lemma1}
Let $D_1,\dots,D_n$ be disks with centers
$O_1,\dots,O_n$ and radii $r_1,\dots,r_n$.
Assume that $O_1,\dots,O_n$ form a convex
polygon in cyclic order. Suppose that
$D_i$ and $D_{i+1}$ intersect with angle
$\Theta_i$, where $1\le i\le n$ and
$D_{n+1}=D_1$ (Figure~\ref{8}).
Assume that
\(
\sum_i(\pi-\Theta_i)>2\pi+\varepsilon_0.
\)
Let $E$ be a bounded set of positive diameter
that meets the interior of the polygon and
intersects every boundary circle $\partial D_i$.
Then there exists a constant
$C=C(\varepsilon_0)>0$ such that
\[
(\operatorname{diam}(E))^{-1/2}
<C\sum_i r_i^{-1/2}.
\]
\end{lemma}

\begin{proof}
Let $\delta=\operatorname{diam}(E)$.
Choose $x\in E$ in the interior of the
convex polygon. For each $i$, choose
$y_i\in E\cap\partial D_i$.
Set
\[
r_i'=d(x,O_i)=r_i+\Delta r_i.
\]
By the reverse triangle inequality,
\[
|\Delta r_i|
=\bigl||x-O_i|-|y_i-O_i|\bigr|
\le |x-y_i|
\le\delta.
\]

If $\delta\ge\frac12\min_i r_i$, then
\[
\delta^{-1/2}
\le\sqrt{2}\sum_i r_i^{-1/2},
\]
and the conclusion follows.
Hence we may assume
\[
\delta<\frac12\min_i r_i.
\]
In particular, $r_i'\ge r_i/2>0$.

Let
\[
\Theta_i'=\pi-\angle O_i xO_{i+1},
\qquad
d_i=|O_i-O_{i+1}|.
\]
By the law of cosines,
\[
\cos\Theta_i
=-\frac{r_i^2+r_{i+1}^2-d_i^2}
{2r_ir_{i+1}},
\]
and
\[
\cos\Theta_i'
=-\frac{r_i'^2+r_{i+1}'^2-d_i^2}
{2r_i'r_{i+1}'}.
\]

Since the polygon is convex and $x$
lies in its interior, we have
\[
\sum_i(\pi-\Theta_i')
=\sum_i\angle O_i xO_{i+1}
=2\pi.
\]

For each $i$, the following inequality holds:
\[
|\Theta_i-\Theta_i'|
\le\frac{\pi}{\sqrt{2}}
\sqrt{|\cos\Theta_i-\cos\Theta_i'|}.
\]

Set
\[
F_i(s,t)=\frac{d_i^2-s^2-t^2}{2st}.
\]
Along the line segment joining
$(r_i,r_{i+1})$ and $(r_i',r_{i+1}')$,
the triangle inequalities hold.
Thus $|F_i(s,t)|\le1$ along this segment.

Moreover,
\[
\frac{\partial F_i}{\partial s}
=-\frac1t-\frac{F_i(s,t)}s,
\qquad
\frac{\partial F_i}{\partial t}
=-\frac1s-\frac{F_i(s,t)}t.
\]
Since $|\Delta r_i|\le\delta$ and
$r_i'\ge r_i/2$, we obtain
\[
|\cos\Theta_i-\cos\Theta_i'|
\le8\delta
\left(\frac1{r_i}+\frac1{r_{i+1}}\right).
\]

Notice that
\[
\begin{aligned}
\sum_i|\Theta_i-\Theta_i'|
&\ge
\sum_i(\pi-\Theta_i)
-\sum_i(\pi-\Theta_i')\\
&>\varepsilon_0.
\end{aligned}
\]

Therefore,
\[
\varepsilon_0
<
2\pi\sqrt{\delta}
\sum_i\sqrt{
\frac1{r_i}+\frac1{r_{i+1}}
}
\le
4\pi\sqrt{\delta}
\sum_i r_i^{-1/2}.
\]
Consequently,
\[
\delta^{-1/2}
<
\frac{4\pi}{\varepsilon_0}
\sum_i r_i^{-1/2}.
\]

Taking
\[
C(\varepsilon_0)
=\max\left\{2,\frac{4\pi}{\varepsilon_0}\right\}
\]
completes the proof.
\end{proof}

\begin{remark}
In this proof, we can actually treat non-intersection edges as intersection with angle $0$. we shall use the actual lengths $r_i'$ for the cosine value estimate.
\end{remark}

\begin{lemma}\label{lemma2}
There exists $\delta_2 = \delta(\varepsilon_1) < 1$ such that:

$$\sum_{v \sim u} r_v - \max_{v \sim u} \{r_v\} > \delta_2 r_u.$$
\end{lemma}

\begin{figure}[htbp]
\centering\includegraphics[scale=0.02]{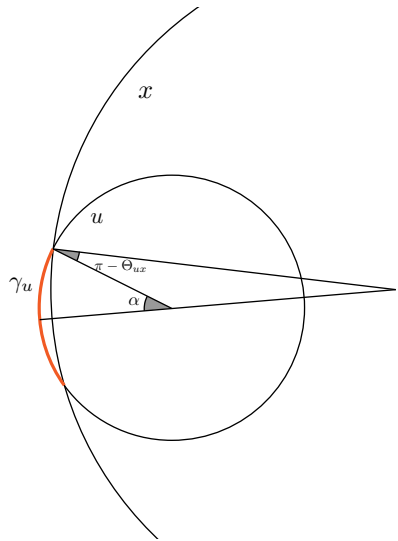}
\caption{The orange arc must be covered}
\label{9}
\end{figure}

\begin{proof}
Let $r_x = \max_{v \sim u} \{r_v\}$. Then the remaining disks cover an arc $\gamma_u$. (Figure~\ref{9}) The length of $\gamma_u$ is $l_u = 2\alpha r_u > 2\varepsilon_1 r_u$. Also, the length $l_u < \sum_{v \sim u, v \neq x} \pi r_v$. Thus, $\sum_{v \sim u} r_v > \frac{2\varepsilon_1}{\pi} r_u$.
\end{proof}

\begin{theorem}\label{thm:re-ring-lemma}
    Let $(G, \Theta)$ be a tame $\IIA$ with $\Theta \in [0, \pi - \varepsilon_1]^E$. Let its embedding on the plane have radius $r: V \rightarrow \mathbb{R}$. Then there exists a uniform constant $C = C(\varepsilon_0, \varepsilon_1)$ such that, for any edge $u \sim v$:

$$\frac{r(v)}{r(u)} > e^{-C\cdot S(u)},$$
where $S(u)=\sum_{v' \sim u} \deg(v')$, which we call it the flower degree of $u$.
\end{theorem}

\begin{lemma}\label{lemma3}
 Let $u \sim a$, $u \sim b$, and $u \sim v_1, \dots, v_n$. Let $a \sim v_1 \sim v_2 \sim \dots \sim v_n \sim b$ be a chain connecting $a$ and $b$. Then there exists $C_3 = C(\varepsilon_0, \varepsilon_1)$ such that:

$$(\sum r_{v_i})^{-1/2} < C_3^{\min(\deg a, \deg b)} (r_u^{-1/2} + r_a^{-1/2} + r_b^{-1/2}).$$
\end{lemma}

\begin{figure}[htbp]
\centering
\includegraphics[scale=0.02]{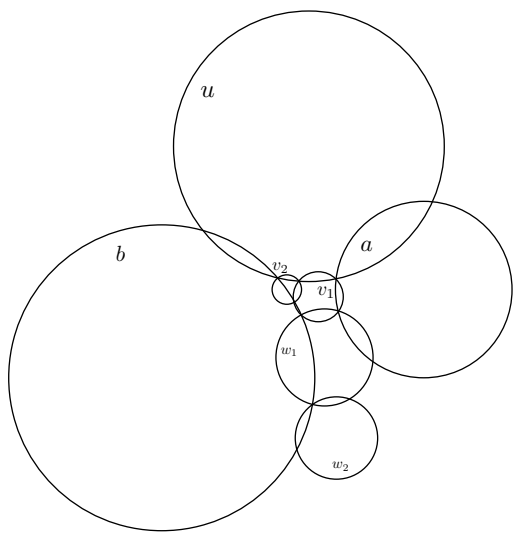}
\caption{A possible arrangement of circles}
\label{10}
\end{figure}

\begin{proof}
Let $E$ be the union of disks $v_1, v_2, \dots, v_n$, and $w_1, w_2, \dots, w_m$ be disks around $b$ (except disks in $E$ and the center $u$, Figure~\ref{10}). Order the 
\(w_i\)'s in blocks around \(b\), each block starting with a disk intersecting \(E\). Since $v_1,\ldots,v_n$ form a chain of pairwise-intersecting
disks, we have $\text{diam}(E) \leq 2\sum_{i=1}^n r_{v_i}$.

If $a$ and $b$ intersect: we are done with $(C_2')$ and $C_3=C(\varepsilon_0)$ from Lemma \ref{lemma1}.

Otherwise, consider disks $w_i$ for $1 \le i \le m$, $w_i \ne a$.

If $w_i \sim a$ : By Lemma \ref{lemma1} and $(C_2')$, $\exists \,C = C(\varepsilon_0) > 1$, 
$$(\text{diam}(E)+2\sum_{j=1}^{i-1}r_{w_j})^{-\frac{1}{2}} < C(r_u^{-\frac{1}{2}} + r_a^{-\frac{1}{2}}+r_b^{-\frac{1}{2}}+r_{w_i}^{-\frac{1}{2}}).$$

If $w_i \nsim a$ : By the remark of Lemma \ref{lemma1} and $\Theta \in [0, \pi - \varepsilon_1]^E$, $\exists C = C(\varepsilon_1) > 1$, 
 
$$(\text{diam}(E)+2\sum_{j=1}^{i-1}r_{w_j})^{-\frac{1}{2}} < C(r_a^{-\frac{1}{2}}+r_b^{-\frac{1}{2}}+r_{w_i}^{-\frac{1}{2}}).$$

We can conclude that there exists $C = C(\varepsilon_0, \varepsilon_1)$, such that 
$$(\text{diam}(E)+2\sum_{j=1}^{i-1}r_{w_j})^{-\frac{1}{2}} < C(r_u^{-\frac{1}{2}} + r_a^{-\frac{1}{2}}+r_b^{-\frac{1}{2}}+r_{w_i}^{-\frac{1}{2}}).$$

Let $C_3 = \frac{8 C^2 + 1}{\delta_2}$, suppose that:

$$(\sum_{i=1}^n r_{v_i})^{-1/2} \ge C_3^{\text{deg } b} (r_u^{-1/2} + r_a^{-1/2} + r_b^{-1/2}).$$

Then:

$$(\sum_i r_{v_i} + r_{w_1})^{-1/2} \ge C_3^{\text{deg } b - 1} (r_u^{-1/2} + r_a^{-1/2} + r_b^{-1/2}).$$

This process continues until
$$(\sum_i r_{v_i} + \sum_j r_{w_j})^{-1/2} \ge C_3^{\text{deg } b - m} (r_u^{-1/2} + r_a^{-1/2} + r_b^{-1/2}).$$

Hence,

$$\sum r_{v_i} + \sum r_{w_j} \le C_3^{-2} r_b < \delta_2 r_b,$$
which contradicts Lemma \ref{lemma2}.

Repeating this logic for $a$, we obtain the following.

$$(\sum r_{v_i})^{-1/2} < C_3^{\min(\text{deg } a, \text{deg } b)} (r_u^{-1/2} + r_a^{-1/2} + r_b^{-1/2})$$

\end{proof}

If $\sum_{i=1}^n r_{v_i} < C_3^{-2 \min(\text{deg } a, \text{deg } b)} r_u/4$, then $\min(r_a, r_b) < 16C_3^{2 \min(\text{deg } a, \text{deg } b)} \sum r_{v_i}$. 
Setting $$C_0 = \frac{16C_3^2 + 1}{\delta_2^2},$$ we can prove by contradiction: 

Suppose that there is $v \sim u$ such that 
$$r_{v} \le C_0^{-\sum_{v'\sim u}\text{deg } v'} r_u,$$
then using Lemma \ref{lemma3}, there exists $v_1 \sim u$, $v_1 \sim v$ that

$$r_{v} + r_{v_1} \le C_0^{-\sum_{v'\sim u, v'\ne v_1}\text{deg } v'} r_u.$$

Again, there exists $v_2 \sim u$, $v_2 \sim v\ \text{or}\ v_1$ such that

$$r_{v} + r_{v_1} + r_{v_2} \le C_0^{-\sum_{v'\sim u, v'\ne v_1,v'\ne v_2}\text{deg } v'} r_u.$$

The process continues until there is only one disk $w$ left:

$$\sum_{v'\sim u, v' \ne w}r_{v'} \le C_0^{-\text{deg } w} r_u\leq C_0^{-2}r_u < \delta_2r_u,$$
which contradicts Lemma~\ref{lemma2}. By taking the constant $C = \log C_0$, we finish the proof.

\begin{remark}
    The figure~\ref{11} shows that the constant $\varepsilon_0$ is crucial in the proof.

\begin{figure}[htbp]
\centering
\includegraphics[scale=0.018]{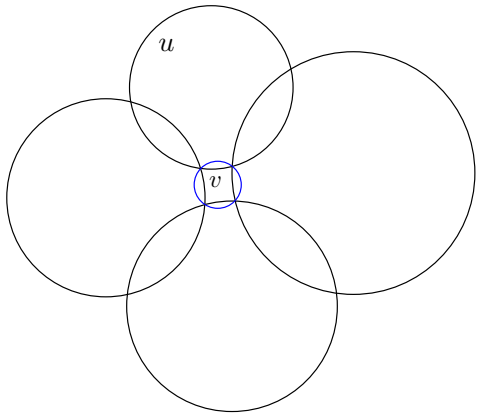}
\caption{The center circle can be arbitrarily small when $\varepsilon_0$ tends to 0}
\label{11}
\end{figure}

\end{remark}

\subsection{Convergence}
The focus of this subsection is to establish the convergence and to analyze the resulting exit distribution, thereby proving items (1) and (2) of Theorem~\ref{thm:icp-poisson}.

\medskip
\begin{lemma}[Exponential decay of radii]\label{lem:ICP-51}
Let $(G,\rho,\Theta)$ be a reversible,
$\ICP$-hyperbolic random rooted tame $\IIA$
with $\sup_e\Theta(e)<\pi$. Assume that
\(
\mathbb E[\deg_G(\rho)^2]<\infty.
\)
Let $(X_n)_{n\ge0}$ be simple random walk
started at $\rho$. Then almost surely
\[
\limsup_{n\to\infty}
\frac{\log r(X_n)}{n}<0.
\]
\end{lemma}

\begin{proof}
We may assume that the reversible law is
ergodic; otherwise, we work on its ergodic
components.

Let $\widehat{\mathbb P}$ denote the
inverse-degree-biased unimodular law:
\[
\frac{d\widehat{\mathbb P}}{d\mathbb P}
=
\frac{\deg_G(\rho)^{-1}}
{\mathbb E[\deg_G(\rho)^{-1}]}.
\]
Write $\widehat{\mathbb E}$ for expectation
under $\widehat{\mathbb P}$. Then
\[
\widehat{\mathbb E}[\deg_G(\rho)^3]
=
\frac{\mathbb E[\deg_G(\rho)^2]}
{\mathbb E[\deg_G(\rho)^{-1}]}
<\infty.
\]
Since $(G,\rho,\Theta)$ is
$\ICP$-hyperbolic, Theorem~\ref{thm:iag-dichotomy},
applied under $\widehat{\mathbb P}$,
yields invariant non-amenability.

\smallskip
By \cite[Theorem~3.2]{AHNR16}, applied
to the inverse-degree-biased unimodular law,
there exists an ergodic percolation $\omega$
such that
\[
\sup_{v\in V(\omega)}\deg_G(v)
\le M<\infty,
\qquad
\mathbf i_E(\omega)>0
\quad\text{a.s.}
\]
We work with a jointly ergodic marked law
including $\Theta$ and $\omega$.

where $\mathbf{i}_E(\omega)$ is the edge Cheeger constant of $\omega$ in the
weighted graph notation with $w\equiv 1$, namely
\[
\mathbf{i}_E(\omega)
=
\inf\Bigl\{
\frac{w(\partial_E W)}{w(W)}
:\ 
\emptyset\neq W\subset V(\omega),\ W\ \text{finite}
\Bigr\}.
\]

Define the visit times to $V(\omega)$ by
\[
N_0
:=
\inf\{n\ge 0:X_n\in V(\omega)\},
\qquad
N_{m+1}
:=
\inf\{n>N_m:X_n\in V(\omega)\},
\quad m\ge 0.
\]
Following \cite{AHNR16}, define a weight function on the induced graph on $V(\omega)$ by
\[
w(u,v)
:=
\deg_G(u)\,
\mathbb{P}_u\!\bigl(X_{N_1}=v\bigr),
\qquad u,v\in V(\omega).
\]
Let $w(u):=\sum_{v}w(u,v)$ be the induced vertex weight, and extend $w(\cdot)$ to sets by
additivity. Then the induced chain
\[
Y_m:=X_{N_m},\qquad m\ge 0,
\]
is a reversible random walk on the weighted graph $(\omega,w)$ with transition probabilities
\[
P_\omega(u,v)
:=
\mathbb{P}_u(Y_1=v)
=
\frac{w(u,v)}{w(u)},
\qquad u,v\in V(\omega),
\]
and reversible measure $\pi_\omega(u):=w(u)$.

\begin{lemma}\label{lem:cheeger-induced}
The edge Cheeger constant of the weighted graph $(\omega,w)$ satisfies
\[
\mathbf{i}_E(\omega,w)\ \ge\ \frac{1}{M}\,\mathbf{i}_E(\omega).
\]
In particular, $\mathbf{i}_E(\omega,w)>0$ almost surely.
\end{lemma}

\begin{proof}
For $u\in V(\omega)$,
\begin{align*}
w(u)
&=\sum_{v\in V(\omega)} w(u,v) \\
&=\sum_{v\in V(\omega)} \deg_G(u)\,\PP_u(Y_1=v) \\
&=\deg_G(u)\sum_{v\in V(\omega)} \PP_u(Y_1=v) \\
&=\deg_G(u)\,\PP_u(N_1<\infty) \\
&\le \deg_G(u) \\
&\le M.
\end{align*}
so $w(u)\le M$ for all $u\in V(\omega)$. Hence for every finite nonempty $W\subset V(\omega)$,
\begin{equation}\label{eq:wW-upper}
w(W)=\sum_{u\in W} w(u)\ \le\ M\,|W|.
\end{equation}

Let $W\subset V(\omega)$ be finite and nonempty, and let $\partial_E W$ be its 
edge boundary in $\omega$. Fix an edge $uv\in\partial_E W$ with $u\in W$ and $v\notin W$.
Since $uv$ is an edge of $\omega$,
\[
\mathbb{P}_u(X_1=v)=\frac{1}{\deg_G(u)}.
\]
On the event $\{X_1=v\}$ the walk hits $V(\omega)$ at time $1$, hence $N_1=1$ and
$Y_1=X_{N_1}=v$. Therefore,
\[
\mathbb{P}_u(Y_1=v)\ \ge\ \mathbb{P}_u(X_1=v)=\frac{1}{\deg_G(u)}.
\]
Multiplying by $\deg_G(u)$ yields
\[
w(u,v)=\deg_G(u)\,\mathbb{P}_u(Y_1=v)\ \ge\ 1.
\]
Summing over all boundary edges gives
\begin{equation}\label{eq:boundary-compare}
w(\partial_E W)
=
\sum_{e\in \partial_E W} w(e)
\ \ge\
|\partial_E W|.
\end{equation}

Combining \eqref{eq:wW-upper} and \eqref{eq:boundary-compare} yields
\[
\frac{w(\partial_E W)}{w(W)}
\ \ge\
\frac{|\partial_E W|}{M\,|W|}.
\]
Since $\omega$ is locally finite and has no isolated vertices on $V(\omega)$,
\[
w_{1}(W):=\sum_{u\in W}\deg_\omega(u)\ \ge\ |W|.
\]
Hence
\[
\frac{|\partial_E W|}{M\,|W|}
\ \ge\
\frac{1}{M}\,
\frac{|\partial_E W|}{\sum_{u\in W}\deg_\omega(u)}
=
\frac{1}{M}\,
\frac{w_{1}(\partial_E W)}{w_{1}(W)}.
\]
Taking the infimum over all finite nonempty $W\subset V(\omega)$ gives
\[
\mathbf{i}_E(\omega,w)
=
\inf_W\frac{w(\partial_E W)}{w(W)}
\ \ge\
\frac{1}{M}\,
\inf_W\frac{w_{1}(\partial_E W)}{w_{1}(W)}
=
\frac{1}{M}\,\mathbf{i}_E(\omega),
\]
as claimed.
\end{proof}

\begin{lemma}\label{lem:sr-induced}
The spectral radius of the induced walk on $(\omega,w)$ is strictly less than $1$.
\end{lemma}

\begin{proof}
The induced chain $P_\omega$ is reversible with respect to $\pi_\omega(u)=w(u)$.
For any finite nonempty $W\subset V(\omega)$, define
\[
Q_\omega(W,W^{\mathrm c})
:=
\sum_{u\in W}\sum_{v\notin W} \pi_\omega(u)\,P_\omega(u,v).
\]
Since $\pi_\omega(u)P_\omega(u,v)=w(u,v)$,
\[
Q_\omega(W,W^{\mathrm c})
=
\sum_{u\in W}\sum_{v\notin W} w(u,v)
=
w(\partial_E W).
\]
Also $\pi_\omega(W)=\sum_{u\in W}\pi_\omega(u)=w(W)$, hence
\[
\frac{Q_\omega(W,W^{\mathrm c})}{\pi_\omega(W)}
=
\frac{w(\partial_E W)}{w(W)}.
\]
By Lemma~\ref{lem:cheeger-induced}, $\mathbf{i}_E(\omega,w)>0$. Therefore the conductance
\[
\Phi_\omega
:=
\inf\Bigl\{
\frac{Q_\omega(W,W^{\mathrm c})}{\pi_\omega(W)}
:\ 
\emptyset\neq W\subset V(\omega),\ W\ \text{finite}
\Bigr\}
=
\mathbf{i}_E(\omega,w)
\]
is strictly positive. Cheeger's inequality for reversible Markov chains
(see Lyons-Peres \cite[Theorem 6.7]{LyPer16})
gives a universal constant $c_0>0$ such that
\[
1-\rho(P_\omega)\ \ge\ c_0\,\Phi_\omega^2.
\]
Since $\Phi_\omega>0$, we obtain $\rho(P_\omega)<1$.
\end{proof}

Fix $c>0$ such that $\rho(P_\omega)\le e^{-c}$. Since $1\le w(u)\le M$ for all $u\in V(\omega)$,
the standard spectral-radius bound for reversible chains yields, for all $x,y\in V(\omega)$ and
$m\ge 1$,
\begin{equation}\label{eq:induced-m-step}
P_\omega^{\,m}(x,y)
\ \le\
\sqrt{\frac{\pi_\omega(y)}{\pi_\omega(x)}}\,\rho(P_\omega)^{\,m}
\ \le\
M^{1/2}e^{-c m}.
\end{equation}
Exactly as in~\cite{AHNR16}, this gives, for every $v\in V(\omega)$ and every $m\ge 1$,
\begin{equation}\label{eq:hit-bound}
\Pr_\rho\big[X_{N_m}=v\big]
\le
M^{1/2}e^{-c m}.
\end{equation}

\smallskip

\begin{lemma}
\label{lem:count-no-angle-gap}
Let $\mathcal P$ be an embedded ideal circle pattern in
$\mathbb D$, with Euclidean centres $z(v)$ and radii $r(v)$.
Suppose that, for some $0<\varepsilon_1<\pi$,
\[
0<\Theta(e)\le \pi-\varepsilon_1
\qquad\text{for every }e\in E.
\]
Then
\[
\sum_{v\in V}r(v)^2\le \frac{2}{\sin^2\varepsilon_1}.
\]
Consequently, for every $\tau>0$,
\[
N(\tau):=\#\{v\in V:r(v)\ge\tau\}
\le \frac{2}{\sin^2\varepsilon_1}\tau^{-2}.
\]
\end{lemma}

\begin{proof}
For each vertex $v$, let
\[
\mathcal M(v):=\bigcup_{e\ni v}\widetilde Q_e
\]
be the union of the embedded quadrilaterals incident to $v$.
We first show that
\[
B\bigl(z(v),r(v)\sin\varepsilon_1\bigr)\subset\mathcal M(v).
\]

The point $z(v)$ lies in the interior of $\mathcal M(v)$.
Every boundary segment of $\mathcal M(v)$ has the form $[p,z(w)]$,
where $w\sim v$ and $p$ is one of the two face points
of $\widetilde Q_{\{v,w\}}$.
In the triangle with vertices $z(v),p,z(w)$,
\[
|z(v)-p|=r(v),
\qquad
\gamma:=\angle z(v)pz(w)
=\pi-\Theta(vw)\ge\varepsilon_1.
\]
If $\gamma\ge\pi/2$, the distance from $z(v)$ to the ray
starting at $p$ through $z(w)$ is $r(v)$.
If $\gamma<\pi/2$, that distance is
$r(v)\sin\gamma\ge r(v)\sin\varepsilon_1$.
In either case,
\[
\operatorname{dist}\bigl(z(v),[p,z(w)]\bigr)
\ge r(v)\sin\varepsilon_1.
\]
Thus the distance from $z(v)$ to $\partial\mathcal M(v)$ is at
least $r(v)\sin\varepsilon_1$, proving the inclusion.

The interiors of the
quadrilaterals $\widetilde Q_e$ are pairwise disjoint.
Each quadrilateral belongs to precisely the two sets $\mathcal M(v)$
corresponding to its endpoints. All these quadrilaterals
are contained in $\mathbb D$. Therefore, by countable
additivity of area,
\[
\begin{aligned}
\pi\sin^2\varepsilon_1\sum_{v\in V}r(v)^2
&\le \sum_{v\in V}\operatorname{area}(\mathcal M(v))\\
&=2\sum_{e\in E}\operatorname{area}(\widetilde Q_e)\\
&\le 2\pi.
\end{aligned}
\]
This gives the first assertion. Finally,
\[
N(\tau)\tau^2
\le \sum_{v:r(v)\ge\tau}r(v)^2
\le \frac{2}{\sin^2\varepsilon_1},
\]
which proves the counting estimate.
\end{proof}

We next convert \eqref{eq:hit-bound} into exponential decay of the radii along the subsequence
$(X_{N_m})$. 
By Lemma~\ref{lem:count-no-angle-gap}, with
$C=2/\sin^2\varepsilon_1<\infty$,
we have
\[
\#\{v:r(v)\ge\tau\}
\le C\tau^{-2},
\qquad 0<\tau<1.
\]

With $\tau=e^{-c m/4}$ this gives
\[
\#\{v:\ r(v)\ge e^{-c m/4}\}\ \le\ C\,e^{c m/2}.
\]
Therefore, using \eqref{eq:hit-bound},
\begin{align*}
\Pr_\rho\bigl[r(X_{N_m})\ge e^{-c m/4}\bigr]
&=
\sum_{v:\, r(v)\ge e^{-c m/4}} \Pr_\rho[X_{N_m}=v] \\
&\le
\#\{v:\ r(v)\ge e^{-c m/4}\}\cdot M^{1/2}e^{-c m} \\
&\le
C e^{c m/2}\cdot M^{1/2}e^{-c m}
=
C M^{1/2}e^{-c m/2}.
\end{align*}
This is summable in $m$, so by Borel--Cantelli, almost surely for all large $m$,
\[
r(X_{N_m})\ \le\ e^{-c m/4}.
\]
Equivalently,
\begin{equation}\label{eq:induced-exp}
\limsup_{m\to\infty}\frac{\log r(X_{N_m})}{m}
\ \le\ 
-\frac{c}{4}
\ <0.
\end{equation}

We now use the refined ring lemma for $\ICP$ to transfer the exponential decay of
radii along the induced walk on the percolation $\omega$ to the original random walk
on $G$.

\begin{lemma}[Integrability]
\label{lem:S-moment}
Let $(G,\rho,\Theta)$ be a unimodular
random rooted IAG satisfying
\begin{equation}
\E\bigl[\deg_G(\rho)^3\bigr]<\infty.
\label{eq:boundary-third-moment}
\end{equation}
Then
\[
    \E[S(\rho)]
    =
    \E\bigl[\deg_G(\rho)^2\bigr]
    <\infty,
\]
and
\[
    \E\bigl[\deg_G(\rho)S(\rho)\bigr]
    \le
    \E\bigl[\deg_G(\rho)^3\bigr]
    <\infty.
\]
\end{lemma}

\begin{proof}
Apply the mass-transport principle to
\[
    F_1(G,u,v)
    :=
    \mathbf 1_{\{u\sim v\}}\deg_G(v).
\]
The total mass sent from \(\rho\) is
\[
    \sum_{v\sim \rho}\deg_G(v)
    =
    S(\rho),
\]
whereas the total mass received at \(\rho\) is
\[
    \sum_{u\sim \rho}\deg_G(\rho)
    =
    \deg_G(\rho)^2.
\]
Hence
\[
    \E[S(\rho)]
    =
    \E\bigl[\deg_G(\rho)^2\bigr]
    <\infty.
\]
For the second assertion, write \(d(v):=\deg_G(v)\). Since
\[
    d(\rho)S(\rho)
    =
    \sum_{v\sim\rho}d(\rho)d(v),
\]
the inequality \(2ab\le a^2+b^2\) gives
\[
\begin{aligned}
    d(\rho)S(\rho)
    &\le
    \frac12
    \sum_{v\sim\rho}
    \bigl(d(\rho)^2+d(v)^2\bigr) \\
    &=
    \frac12 d(\rho)^3
    +
    \frac12\sum_{v\sim\rho}d(v)^2.
\end{aligned}
\]
Now apply the mass-transport principle to
\[
    F_2(G,u,v)
    :=
    \mathbf 1_{\{u\sim v\}}\deg_G(v)^2.
\]
The total mass sent from \(\rho\) is
\[
    \sum_{v\sim\rho}\deg_G(v)^2,
\]
whereas the total mass received at \(\rho\) is
\[
    \sum_{u\sim\rho}\deg_G(\rho)^2
    =
    \deg_G(\rho)^3.
\]
Therefore,
\[
    \E\left[\sum_{v\sim\rho}\deg_G(v)^2\right]
    =
    \E\bigl[\deg_G(\rho)^3\bigr].
\]
Combining the preceding estimates yields
\[
    \E\bigl[\deg_G(\rho)S(\rho)\bigr]
    \le
    \E\bigl[\deg_G(\rho)^3\bigr]
    <\infty.
\]
\end{proof}

By stationarity and ergodicity of the joint
marked environment $(G,\rho,\Theta,\omega)$,
the density of visits to $V(\omega)$ is
\[
q:=\mathbb P(\rho\in V(\omega))>0.
\]
Indeed, the ergodic theorem gives
\[
\frac1n\sum_{i=0}^{n-1}
\mathbf 1_{\{X_i\in V(\omega)\}}
\longrightarrow q
\qquad\text{a.s.}
\]

Consequently,
\begin{equation}\label{eq:linear-Nm-ICP-final}
\frac{N_m}{m}\longrightarrow
q^{-1}\in[1,\infty)
\qquad\text{a.s.}
\end{equation}
In particular, setting
$m(n):=\max\{m:N_m\le n\}$, we obtain
\begin{equation}\label{eq:Nm-over-n-ICP-final}
\frac{N_{m(n)}}{n}\longrightarrow1
\qquad\text{a.s.}
\end{equation}

From \eqref{eq:induced-exp} and \eqref{eq:linear-Nm-ICP-final} we obtain
\begin{equation}\label{eq:limsup-Nm-ICP-final}
\limsup_{m\to\infty}\frac{\log r(X_{N_m})}{N_m}
=
\limsup_{m\to\infty}\Bigl(\frac{m}{N_m}\Bigr)\Bigl(\frac{\log r(X_{N_m})}{m}\Bigr)
\ \le\ -\frac{cq}{4}\ <\ 0.
\end{equation}
Combining \eqref{eq:limsup-Nm-ICP-final} with \eqref{eq:Nm-over-n-ICP-final} yields
\begin{equation}\label{eq:subseq-neg-ICP-final}
\limsup_{n\to\infty}\frac{\log r(X_{N_{m(n)}})}{n}\ <\ 0.
\end{equation}

By the refined Ring Lemma, for every edge $u\sim v$,
\[
\log\frac{r(v)}{r(u)}\le A S(v).
\]
Consequently,
\begin{equation}\label{eq:path-sum-ICP-final}
  \log\frac{r(X_n)}{r(X_{N_{m(n)}})}
\le A\sum_{i=N_{m(n)}+1}^{n}S(X_i).
\end{equation}

Hence
\begin{equation}\label{eq:main-split-ICP-final}
  \frac{\log r(X_n)}{n}
\le
\frac{\log r(X_{N_{m(n)}})}{n}
+
\frac{A}{n}\sum_{i=N_{m(n)}+1}^{n}S(X_i).  
\end{equation}

It remains to show that the second term in
\eqref{eq:main-split-ICP-final} vanishes.

By reversibility and ergodicity, the
environment seen from the walk
$(G,X_n,\Theta)$ is stationary and ergodic
under the annealed law.

Moreover, by reversibility and the
Cauchy--Schwarz inequality,
\[
\begin{aligned}
\mu:=\mathbb E[S(\rho)]
&=\mathbb E[
\deg_G(\rho)\deg_G(X_1)]\\
&\le
\sqrt{
\mathbb E[\deg_G(\rho)^2]\,
\mathbb E[\deg_G(X_1)^2]}\\
&=\mathbb E[\deg_G(\rho)^2]
<\infty.
\end{aligned}
\]

Therefore, Birkhoff's ergodic theorem gives
\begin{equation}\label{eq:birkhoff-S}
\frac{1}{n}\sum_{i=0}^{n-1}S(X_i)
\xrightarrow[n\to\infty]{\mathrm{a.s.}}
\mu.
\end{equation}

Set $k_n:=N_{m(n)}$.
By \eqref{eq:Nm-over-n-ICP-final},
\[
k_n\to\infty,
\qquad
\frac{k_n}{n}\to1
\quad\text{almost surely}.
\]
Now
\[
\begin{aligned}
\frac{A}{n}\sum_{i=k_n+1}^{n}S(X_i)
={}&
A\frac{n+1}{n}
\left(\frac{1}{n+1}\sum_{i=0}^{n}S(X_i)\right)\\
&-
A\frac{k_n+1}{n}
\left(\frac{1}{k_n+1}
\sum_{i=0}^{k_n}S(X_i)\right).
\end{aligned}
\]

By \eqref{eq:birkhoff-S}, both averages in
parentheses converge almost surely to $\mu$,
the second along the random subsequence $k_n\to\infty$.
Moreover,
\(
\frac{k_n+1}{n}\to1.
\)
Consequently,
\begin{equation}\label{eq:error-vanish-ICP-final}
\frac{A}{n}\sum_{i=N_{m(n)}+1}^{n}
S(X_i)
\xrightarrow[n\to\infty]{\mathrm{a.s.}}0.
\end{equation}

Combining \eqref{eq:main-split-ICP-final},
\eqref{eq:subseq-neg-ICP-final}, and
\eqref{eq:error-vanish-ICP-final}, we obtain
\[
\limsup_{n\to\infty}
\frac{\log r(X_n)}{n}<0
\qquad\text{almost surely}.
\]
\end{proof}

We now prove Theorem~\ref{thm:icp-poisson}.
By degree biasing, it suffices to work
under the equivalent reversible law
\[
\frac{d\mathbb P^{\mathrm{rev}}}{d\mathbb P}
=
\frac{\deg_G(\rho)}
{\mathbb E[\deg_G(\rho)]}.
\]
Indeed,
\[
\mathbb E^{\mathrm{rev}}[\deg_G(\rho)^2]
=
\frac{\mathbb E[\deg_G(\rho)^3]}
{\mathbb E[\deg_G(\rho)]}
<\infty.
\]
Thus Lemma~\ref{lem:ICP-51} applies under
$\mathbb P^{\mathrm{rev}}$.
We henceforth work under the reversible law
for stationary random-walk arguments, and
under the original unimodular law when
applying the Mass-Transport Principle.
Since the two laws are equivalent,
almost-sure statements transfer between them.

\begin{proposition}
Almost surely the Euclidean centers
$z(X_n)$ and hyperbolic centers $z_h(X_n)$
converge to a common boundary point
$\delta\in\partial\D$
\end{proposition}

\begin{proof}
For the SRW $(X_i)$ on $G$, let $z(v)$ be the Euclidean centre and $r(v)$ the
radius of the circle corresponding to $v$ in $\D$.
The polyline length along successive centres satisfies
\begin{align*}
L
\;:=\;
\sum_{i\ge0}\bigl| z(X_{i+1}) - z(X_i) \bigr|
\;\le\;
\sum_{i\ge0}\bigl( r(X_i) + r(X_{i+1}) \bigr).
\end{align*}
Hence
\begin{align*}
L
&\le
r(\rho)
\;+\;
2\sum_{i\ge1} r(X_i)
\;+\;
\lim_{n\to\infty} r(X_n).
\end{align*}
Since in the hyperbolic case $r(X_n)\to0$ a.s., we obtain
\[
L
\;\le\;
r(\rho)
\;+\;
2\sum_{i\ge1} r(X_i),
\]
which is almost surely finite by Lemma~\ref{lem:ICP-51}.
Hence $(z(X_n))$ is Cauchy and converges to some point in $\overline{\D}$.
Since $G$ is VEL-hyperbolic, simple
random walk on $G$ is transient.
Moreover, the ICP is locally finite
in $\mathbb D$, so every compact subset
of $\mathbb D$ meets only finitely many
circles. Consequently, the walk cannot
have its circle centres converging to
an interior point of $\mathbb D$.
Hence
\[
z(X_n)\ \longrightarrow\ \delta\in\partial\D
\qquad\text{a.s.}
\]
Moreover, for every $v\in V(G)$ the hyperbolic centre lies in the Euclidean disc
$B\!\bigl(z(v),r(v)\bigr)$, hence
\[
\bigl|z_h(v)-z(v)\bigr|\ \le\ r(v).
\]
Since $r(X_n)\to0$ and $z(X_n)\to\delta\in\partial\D$ a.s., we also have
\[
z_h(X_n)\ \longrightarrow\ \delta
\qquad\text{a.s.}
\]
\end{proof}

\subsection{Non-atomicity and full support of exit measure}

Fix $u\in V(G)$. By ICP rigidity,
the realization in $\D$ is unique up
to a M\"obius automorphism of $\D$.
Normalize the realization by sending
the hyperbolic centre $z_h(u)$ to $0$.
Any two such normalized realizations
differ by a rotation of $\D$ (or a
reflection if orientation-reversing
isometries are allowed). Thus the
Euclidean distances between circle
centres are independent of the
remaining choice.
Denote the Euclidean centre map in
this normalization by $z^{(u)}$,
and define a metric on $V(G)$ by
\[
d_G^{(u)}(v,w)
\;:=\;
\bigl|\,z^{(u)}(v)-z^{(u)}(w)\,\bigr|.
\]
If we change the basepoint from $u$ to $v$, the two normalized patterns differ by
a disc automorphism; hence their boundaries are canonically homeomorphic.
Thus $\{d_G^{(u)}\}_{u\in V}$ is a compatible family of metrics in the sense
of \cite[\S5.2]{AHNR16}.  Let $\overline V=V\sqcup\partial V$ be the resulting
completion and boundary.

Since $z^{(u)}(X_n)\to\delta\in\partial\D$ a.s., the sequence $(X_n)$ is
$d_G^{(u)}$-Cauchy and hence converges a.s.\ in the completion $\overline V$.
Moreover, under the canonical embedding induced by $z^{(u)}$, the limit point of
$X_n$ in $\partial V$ is mapped to $\delta\in\partial\D$. The following argument of Lemma \ref{lem:atom-or-none}, \ref{lem:hit-past-ICP} and \ref{lem:hit-finite} are adapted from the method of Angel–Hutchcroft–Nachmias–Ray \cite{AHNR16}.

\begin{lemma}
\label{lem:atom-or-none}
Let $d=\{d_G^{(u)}\}$ be any compatible family of metrics, and let $(G,\rho,\Theta)$ be a
stationary random rooted $\IIA$ graph.
If $X_n$ converges a.s.\ to $\partial V$ in the completion defined by $d$, then the
exit measure on $\partial V$ is almost surely either trivial (a single atom of mass $1$)
or non-atomic.
\end{lemma}

\begin{proof}
Fix the $(G,\rho,\Theta)$. For each $\xi\in\bdry V$ define the harmonic function
\[
h_\xi(v)
\;:=\;
\PP_v\!\left(\lim_{n\to\infty}X_n=\xi\right).
\]
By L\'evy's $0$--$1$ law, along the walk we have
\[
h_\xi(X_n)\ \longrightarrow\ \1\{\lim_{n\to\infty}X_n=\xi\}
\qquad\text{a.s.}
\]
Let
\[
M(G,v)
\;:=\;
\sup_{\xi\in\bdry V} h_\xi(v)
\]
denote the maximal atom mass. Because the boundary topology of the invariant
completion is independent of the root, the process $M(G,X_n)$ is stationary.

Let $\mathcal A$ be the set of atoms of
the exit measure from $\rho$.
Since the walk is irreducible, exit
measures from different vertices
are mutually absolutely continuous;
thus $\mathcal A$ is independent of
the starting vertex. Moreover,
$\mathcal A$ is countable.

Define
\[
q(v):=\sum_{\xi\in\mathcal A}h_\xi(v).
\]
By L\'evy's $0$--$1$ law,
\[
q(X_n)\longrightarrow
\mathbf 1_{\{\lim X_n\in\mathcal A\}}
\qquad\text{a.s.}
\]
On the event $\{\lim X_n\in\mathcal A\}$,
we have
\[
h_{\lim X_n}(X_n)\longrightarrow1
\qquad\text{a.s.},
\]
by applying L\'evy's law
simultaneously to the countable
family of atoms. Therefore
$M(G,X_n)\to1$ on this event.
On its complement,
\[
0\le M(G,X_n)\le q(X_n)\to0.
\]
Consequently,
\[
M(G,X_n)\longrightarrow
\mathbf 1_{\{\lim X_n\in\mathcal A\}}
\qquad\text{a.s.}
\]
Therefore
\[
M(G,X_n)\ \longrightarrow\ \1\{\lim X_n\text{ is an atom}\}
\qquad\text{a.s.}
\]
Stationarity then implies $M(G,\rho)\in\{0,1\}$ almost surely. Hence either there
are no atoms, or there is a single atom of mass $1$.
\end{proof}

Let $\mu=\mu_{G,\rho}$ denote the exit measure on $\partial V$, i.e.\ the conditional
law of $\lim_{n\to\infty}X_n$ given $(G,\rho,\Theta)$.
Applying Lemma~\ref{lem:atom-or-none} to the $\ICP$-compatible family $\{d_G^{(u)}\}$
shows that $\mu$ is almost surely either trivial or non-atomic.
It remains to rule out the trivial case.

Assume for contradiction that $\mu$ is almost surely a single atom $\{\xi\}$ of mass $1$.
Choose a M\"obius map $\Psi:\D\to\Hh=\{w\in\C:\Im w>0\}$ with $\Psi(\xi)=\infty$, and push
the pattern forward to $\Psi(\mathcal{P})$ in $\Hh$.
In this normalization the boundary point $\infty$ is fixed, so the remaining freedom is
\[
w\ \longmapsto\ a w + b,
\qquad a>0,\ b\in\mathbb{R}.
\]
In particular, the straight-line embedding obtained by joining Euclidean circle centres
is well-defined up to similarities. By Theorem \ref{thm:section3-1} we have
\begin{equation}\label{eq:angle-expectation}
\EE\bigl[T(\rho)\bigr]
\;=\;
\EE\bigl[\kappa(\rho)\bigr].
\end{equation}
On the other hand, in the straight-line embedding of $\Psi(\mathcal{P})\subset\Hh$, every vertex has
total angle $2\pi$. Transporting angles and applying unimodularity yields
\begin{equation}\label{eq:2pi}
\EE\bigl[\kappa(\rho)\bigr]
\;=\;
0.
\end{equation}
Combining \eqref{eq:angle-expectation}--\eqref{eq:2pi} gives $\EE[T(\rho)]=0$.

But the $\ICP$ dichotomy (Theorem \ref{thm:iag-dichotomy}) asserts that in the $\ICP$-hyperbolic case
\[
\EE[T(\rho)]
\;<\;
0,
\]
while equality corresponds to the parabolic case. This contradiction shows that the exit
measure cannot be a single atom. Hence $\mu$ is almost surely non-atomic.

It remains to prove that $\supp(\mu)=\partial\D$ almost surely.
Assume for contradiction that $\supp(\mu)\subsetneq\partial\D$.
Then $\partial\D\setminus\supp(\mu)$ is a non-empty open subset of $\partial\D$,
hence it can be written as a countable disjoint union of open arcs
\[
\partial\D\setminus\supp(\mu)
\;=\;
\bigcup_{i\in I}(\theta_i,\psi_i),
\qquad I\ \text{countable},
\]
and $\mu(\{e^{i\theta_i}\})=\mu(\{e^{i\psi_i}\})=0$ for every $i\in I$ since $\mu$ is non-atomic.

Fix one arc $(\theta,\psi):=(\theta_i,\psi_i)$.
Let $\gamma$ be the hyperbolic geodesic in $\D$ with endpoints
$e^{i\theta}$ and $e^{i\psi}$, and let $S(\theta,\psi)$ be the closed region bounded by
$\gamma$ and the boundary arc $\{e^{it}:t\in[\theta,\psi]\}$.
Set
\[
A
\;:=\;
\bigl\{\,u\in V(G):\ z_h(u)\in S(\theta,\psi)\,\bigr\}.
\]

For each $u\in A$, let $\gamma_u$ be the hyperbolic geodesic ray starting at $z_h(u)$
and heading to $e^{i\theta}$, parameterized by hyperbolic arclength $\gamma_u(t)$, $t\ge 0$.
Define $\mathsf{rec}(u)$ to be the first vertex $v$ such that the circle of $v$
intersects both $\gamma$ and $\gamma_u$, where ``first'' means minimizing
\[
t(v)
\;:=\;
\inf\bigl\{\,t\ge 0:\ \text{the circle of $v$ intersects $\gamma_u(t)$}\,\bigr\}.
\]
In case of ties, break them by a fixed measurable rule (e.g.\ by assigning i.i.d.\ labels
to vertices and choosing the smallest label among minimizers); if no such $v$ exists,
set $\mathsf{rec}(u)=\varnothing$.  To ensure isomorphism-invariance, define \(\operatorname{rec}(u)\) simultaneously for all connected components \(I_j=(\theta_j,\psi_j)\) of \(\partial\mathbb D\setminus\operatorname{supp}(\mu)\). For each \(j\), let \(\gamma_j\) be the hyperbolic geodesic in \(\mathbb D\) with endpoints \(e^{i\theta_j}\) and \(e^{i\psi_j}\), and let \(S_j\) be the closed region bounded by \(\gamma_j\) and the boundary arc \(I_j\). If \(z_h(u)\) lies in the interior of the corresponding region \(S_j\), define \(\operatorname{rec}(u)\) by the preceding rule, using the geodesic \(\gamma_j\) and the ray from \(z_h(u)\) to \(e^{i\theta_j}\). Otherwise, set \(\operatorname{rec}(u)=\varnothing\). Define
\[
m(u,v):=\mathbf 1_{\{\operatorname{rec}(u)=v\}}.
\]
The fixed component \(I_j=(\theta,\psi)\) is used only in the subsequent argument to exhibit a vertex receiving infinite mass. Then each vertex sends at most one unit of mass:
\[
\sum_{v} m(u,v)\ \le\ 1
\qquad\text{for all }u\in V(G).
\]
Obviously, $(G,\Theta,u,v)\mapsto m(u,v)$ is measurable and isomorphism-invariant, hence admissible
for the mass-transport principle.
\begin{figure}[htbp]
\centering
\includegraphics[scale=0.023]{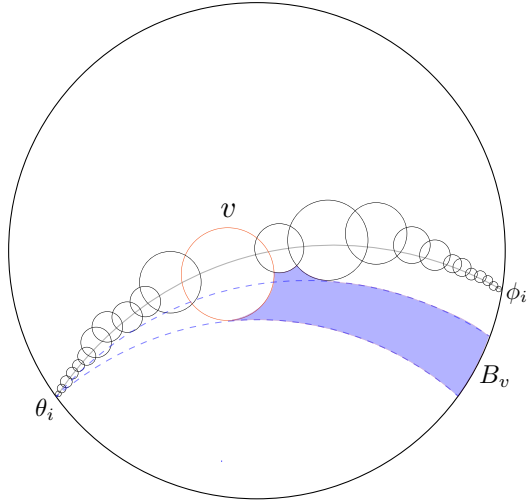}
\caption{The vertex v receives mass from circles with hyperbolic centres
in the shaded area}
\label{12}
\end{figure}
For each $\phi\in(\theta,\psi)$, let $\gamma_\phi$ denote the hyperbolic geodesic
from $e^{i\phi}$ to $e^{i\theta}$.
For a vertex $v$ whose circle intersects $\gamma$, define (Figure~\ref{12})
\[
B_v
\;:=\;
\Bigl\{\,
\phi\in(\theta,\psi):
\ \text{$v$ is the first circle met by $\gamma_\phi$
that also intersects $\gamma$}
\,\Bigr\}.
\]
Up to a null set of exceptional $\phi$ where ties occur, the sets $\{B_v\}$ form a Borel
partition of $(\theta,\psi)$, and hence
\[
(\theta,\psi)
\;=\;
\bigcup_{v:\,\mathrm{circle}(v)\cap\gamma\neq\varnothing} B_v
\qquad\text{(disjoint up to null sets).}
\]
Since there are only countably many vertices, there exists a vertex $v_0$ with
\[
|B_{v_0}|\ >\ 0,
\]
where $|\cdot|$ denotes Lebesgue measure on $(\theta,\psi)$.
Fix $\phi\in B_{v_0}$ outside the exceptional null set.
Then, by definition of $B_{v_0}$, the geodesic $\gamma_\phi$ hits the circle of $v_0$
before hitting any other circle that intersects $\gamma$.
Consequently, every vertex $u\in A$ whose hyperbolic centre lies sufficiently close
to the boundary inside the sector determined by $\gamma_\phi$ and the arc $(\theta,\psi)$
satisfies $\mathsf{rec}(u)=v_0$.
Since circles accumulate at the boundary inside this sector, there are infinitely many
such vertices $u$, and therefore
\[
\sum_{u} m(u,v_0)\ =\ \infty.
\]

This contradicts the (mass-transport principle) unimodularity, since the construction
gives an admissible transport in which each vertex sends at most one unit of mass, but
some vertex receives infinitely many units.
Hence $\partial\D\setminus\supp(\mu)=\varnothing$, i.e.
\[
\supp(\mu)\;=\;\partial\D
\qquad\text{almost surely.}
\]

\subsection{Poisson boundary}

Two results have already been established earlier in this section:
\begin{itemize}
\item[(A)] \textbf{Convergence and exit measure.}
Almost surely $z(X_n)\to \delta^+\in\partial\D$ and $z(X_{-n})\to \delta^-\in\partial\D$.
Moreover the exit measure on $\partial\D$ is non-atomic and has full support.

\item[(B)] \textbf{Exponential decay of radii.}
Almost surely
\[
\limsup_{n\to\infty}\frac{1}{n}\log r(X_n)\ <\ 0.
\]
\end{itemize}

\medskip

Let $\Omega:=V(G)^{\mathbb N_0}$ denote the space of one-sided infinite paths in $G$.
Write $\theta(x_0,x_1,\ldots)=(x_1,x_2,\ldots)$ for the shift, and let
\[
\mathcal I \ :=\ \{A\subset\Omega:\ \theta^{-1}A=A\}
\]
be the shift-invariant $\sigma$-algebra.
As in \cite[§5.3]{AHNR16}, there is an $L^\infty$--isomorphism between bounded harmonic
functions $h$ on $G$ and $L^\infty(\Omega,\mathcal I)$ via
\[
h\ \longleftrightarrow\ g,\qquad
g(x_0,x_1,\ldots)=\lim_{n\to\infty} h(x_n),\qquad
h(v)=\E^v\!\big[g(X_0,X_1,\ldots)\big],
\]
where the limit exists by the bounded martingale convergence theorem; the fact that the two
maps are inverses follows from martingale convergence together with Lévy's $0$--$1$ law.

\begin{proposition}\label{prop:poisson}
Let $G$ be a locally finite graph and let $(X_n)_{n\ge0}$ be simple random walk on $G$.
Assume there exists a boundary random variable $\delta^+\in\partial\D$ such that
$z(X_n)\to\delta^+$ almost surely.
Let $\Omega:=V(G)^{\mathbb N_0}$, let $\theta$ be the left shift, and let $\mathcal I$
be the shift-invariant $\sigma$-algebra. Then the following are equivalent:
\begin{enumerate}
\item[\textnormal{(i)}] $\partial\D$ realises the Poisson boundary, i.e.\ every bounded harmonic
$h$ on $G$ admits a representation
\[
h(v)\;=\;\E^v\!\big[\tilde g(\delta^+)\big]
\qquad\text{for some bounded Borel }\tilde g:\partial\D\to\mathbb R.
\]
\item[\textnormal{(ii)}] For every $A\in\mathcal I$ there exists a Borel set $B\subset\partial\D$
such that
\[
\Pp^\rho\!\big(A\ \triangle\ \{\delta^+\in B\}\big)\;=\;0,
\]
where $\triangle$ denotes symmetric difference.
\end{enumerate}
\end{proposition}

\begin{proof}
For bounded harmonic $h$, the process $(h(X_n))_{n\ge0}$ is a bounded martingale, hence
\[
g(\omega)\;:=\;\lim_{n\to\infty} h(X_n(\omega))
\]
exists almost surely. Moreover $g\circ\theta=g$ a.s., so $g\in L^\infty(\Omega,\mathcal I)$.

\smallskip

Conversely, if $g\in L^\infty(\Omega,\mathcal I)$, define
\[
h(v)\;:=\;\E^v\!\big[g(X_0,X_1,\ldots)\big].
\]
Then $h$ is bounded and harmonic. Writing $\mathcal F_n:=\sigma(X_0,\ldots,X_n)$,
the Markov property gives $h(X_n)=\E^\rho[g\,|\,\mathcal F_n]$, and the martingale convergence theorem
together with Lévy's $0$--$1$ law yields
\begin{equation}\label{eq:Levy}
h(X_n)\ =\ \E^\rho[g\,|\,\mathcal F_n]
\ \xrightarrow[n\to\infty]{\text{a.s.}}\ g.
\end{equation}
Thus bounded harmonic functions are isometrically isomorphic to $L^\infty(\mathcal I)$ via $h\leftrightarrow g$.

\medskip
\noindent(i) $\Rightarrow$ (ii).
Fix $A\in\mathcal I$ and set $g:=\mathbf 1_A$.
By (i), there exists a bounded Borel $\tilde g:\partial\D\to[0,1]$ such that
\[
h(v)\;:=\;\Pp^v(A)\;=\;\E^v\!\big[\tilde g(\delta^+)\big].
\]
On the one hand, \eqref{eq:Levy} (with $g=\mathbf 1_A$) gives $h(X_n)\to \mathbf 1_A$ a.s.
On the other hand, $h(X_n)=\E^\rho[\tilde g(\delta^+)\mid\mathcal F_n]$ is a bounded martingale, hence
$h(X_n)\to \tilde g(\delta^+)$ a.s. Therefore $\mathbf 1_A=\tilde g(\delta^+)$ a.s.
Modifying $\tilde g$ on a $\Pp^\rho$-null set, we may take $\tilde g$ to be $\{0,1\}$-valued.
Let $B:=\{\xi\in\partial\D:\ \tilde g(\xi)=1\}$.
Then $\Pp^\rho\big(A\triangle\{\delta^+\in B\}\big)=0$.

\medskip
\noindent(ii) $\Rightarrow$ (i).
Let $g\in L^\infty(\mathcal I)$ and choose simple functions
\[
g_m=\sum_{j=1}^{J_m} a_{m,j}\,\mathbf 1_{A_{m,j}},
\qquad A_{m,j}\in\mathcal I,
\]
such that $\|g_m-g\|_\infty\to 0$.
For each $(m,j)$, apply (ii) to obtain a Borel set $B_{m,j}\subset\partial\D$ with
\[
\Pp^\rho\big(A_{m,j}\triangle\{\delta^+\in B_{m,j}\}\big)=0,
\]
and define
\[
\tilde g_m(\xi)\;:=\;\sum_{j=1}^{J_m} a_{m,j}\,\mathbf 1_{B_{m,j}}(\xi),
\qquad \xi\in\partial\D.
\]
Then $g_m=\tilde g_m(\delta^+)$ $\Pp^\rho$-a.s. In particular,
\[
\|\tilde g_m(\delta^+)-\tilde g_\ell(\delta^+)\|_{L^\infty(\Pp^\rho)}
=\|g_m-g_\ell\|_\infty,
\]
so $(\tilde g_m(\delta^+))$ is Cauchy in $L^\infty(\Pp^\rho)$ and hence converges (along a subsequence,
and thus along the whole sequence after modification on a null set) to some bounded random variable
$\tilde g(\delta^+)$ that is measurable with respect to $\sigma(\delta^+)$.
Therefore there exists a bounded Borel function $\tilde g:\partial\D\to\mathbb R$ such that
\[
\tilde g(\delta^+)=\lim_{m\to\infty}\tilde g_m(\delta^+)
\qquad\text{$\Pp^\rho$-a.s.}
\]

\smallskip

For every $v$, define $h(v):=\E^v[g]$ and $h_m(v):=\E^v[g_m]$.
Then $\|h_m-h\|_\infty\le \|g_m-g\|_\infty\to 0$.
Since $|\tilde g_m|\le \|g\|_\infty$, dominated convergence gives
\[
h(v)
=\lim_{m\to\infty} h_m(v)
=\lim_{m\to\infty}\E^v\!\big[\tilde g_m(\delta^+)\big]
=\E^v\!\big[\tilde g(\delta^+)\big].
\]
Thus $\partial\D$ realises the Poisson boundary.
\end{proof}

\medskip

We also use the following consequence of optional stopping \cite[Eq.\ (5.5)]{AHNR16}.
If $h\ge 0$ is bounded harmonic and $W\subset V$ with
\[
T_W=\inf\{n\ge 0:\ X_n\in W\},
\]
then
\begin{equation}\label{eq:OS-ineq}
h(v)\ \ge\ \E^v\!\big[h(X_{T_W})\,\mathbf 1_{\{T_W<\infty\}}\big]
\ \ge\ \Pp^v(\mathrm{hit}\ W)\cdot \inf_{u\in W} h(u),
\qquad v\in V.
\end{equation}
Indeed, $h(X_{n\wedge T_W})$ is a bounded martingale, hence
\[
h(v)=\E^v[h(X_{n\wedge T_W})]
\ \ge\ \E^v\!\big[h(X_{T_W})\,\mathbf 1_{\{T_W\le n\}}\big].
\]
Letting $n\to\infty$ and using $h\ge 0$ gives \eqref{eq:OS-ineq}.

\begin{lemma}
\label{lem:hit-past-ICP}
Let $(X_n)_{n\in\Z}$ be the reversible bi-infinite SRW on $G$.
Almost surely,
\[
\Pp^{X_n}\!\big(\mathrm{hit}\ \{X_{-1},X_{-2},\ldots\}\big)\ \xrightarrow[n\to\infty]{}\ 0.
\]
\end{lemma}

\begin{proof}
By (A), the limits $\delta^+=\lim_{n\to\infty} z(X_n)$ and
$\delta^-=\lim_{n\to\infty} z(X_{-n})$ exist, and non-atomicity implies
$\delta^+\neq \delta^-$ almost surely.

\smallskip

Let $\{U_i\}_{i\in I}$ be a countable basis of open arcs in $\partial\D$ with rational endpoints.
For each $i$, define
\[
h_i(v)\ :=\ \Pp^v\!\big(\delta^+\in U_i\big),
\qquad v\in V,
\]
which is bounded harmonic. By Lévy's $0$--$1$ law,
\[
h_i(X_n)\ \longrightarrow\ \mathbf 1_{\{\delta^+\in U_i\}}
\qquad\text{a.s., for all $i$.}
\]

\smallskip

Choose $i_0$ so that $\delta^-\in U_{i_0}$ but $\delta^+\notin U_{i_0}$, and set $h:=h_{i_0}$.
Then almost surely $h(X_n)\to 0$.
Since $\delta^-\in U_{i_0}$, applying Lévy's $0$--$1$ law to the time-reversed walk
$(X_{-n})_{n\ge 0}$ yields $h(X_{-m})\to 1$ a.s., hence $h(X_{-m})\ge \tfrac12$ for all large $m$.
For each fixed $m$ we also have $h(X_{-m})>0$ a.s. because
\[
h(v)=\Pp^v(\delta^+\in U_{i_0})
\ \ge\ \Pp^v(\mathrm{hit}\ \rho)\,\Pp^\rho(\delta^+\in U_{i_0}),
\]
and $\Pp^v(\mathrm{hit}\ \rho)>0$ for every $v$ while $\Pp^\rho(\delta^+\in U_{i_0})>0$ by full support.
Therefore
\[
a:=\inf_{m\ge 1} h(X_{-m})\ >\ 0
\qquad\text{a.s.}
\]
Apply \eqref{eq:OS-ineq} with $W=\{X_{-1},X_{-2},\ldots\}$ and $v=X_n$:
\[
h(X_n)\ \ge\ a\cdot \Pp^{X_n}\!\big(\mathrm{hit}\ \{X_{-1},X_{-2},\ldots\}\big).
\]
Since $h(X_n)\to 0$, the displayed probability tends to $0$.
\end{proof}

\begin{lemma}\label{lem:hit-finite}
Let $(X_n)_{n\in\Z}$ be the reversible bi-infinite SRW on $G$.
Almost surely, for every finite $F\subset V$,
\[
\Pp^{X_n}\!\big(\mathrm{hit}\ F\big)\ \xrightarrow[n\to\infty]{}\ 0.
\]
\end{lemma}

\begin{proof}
Fix a finite $F\subset V$. Choose an open arc $U\subset\partial\D$ with
$\delta^-\in U$ and $\delta^+\notin U$, and define
\[
h(v):=\Pp^v(\delta^+\in U).
\]
Then $h$ is bounded harmonic and $h(X_n)\to 0$ a.s. by Lévy's $0$--$1$ law.
Since $\Pp^\rho(\delta^+\in U)>0$ by full support and $\Pp^u(\mathrm{hit}\ \rho)>0$ for every $u$,
we have $h(u)>0$ for all $u\in F$ and hence
\[
a:=\min_{u\in F} h(u)\ >\ 0.
\]
Applying \eqref{eq:OS-ineq} with $W=F$ and $v=X_n$ gives
\[
h(X_n)\ \ge\ a\cdot \Pp^{X_n}(\mathrm{hit}\ F),
\]
so $\Pp^{X_n}(\mathrm{hit}\ F)\to 0$.
\end{proof}

\medskip

Let $A$ be a shift-invariant event with $\Pp^\rho(A)>0$ and define
\[
h(v)\ :=\ \Pp^v(A).
\]
We will build a Borel set $B\subset\partial\D$ such that
\[
\Pp^\rho\big(A\ \triangle\ \{\delta^+\in B\}\big)\ =\ 0.
\]
By (B), on an almost sure event there exists $c_0>0$ such that
\[
r(X_n)\ \le\ e^{-c_0 n}
\qquad\text{for all sufficiently large $n$.}
\]

\smallskip


Let $\mathcal W_\rho$ be the space of infinite nearest-neighbour paths
$\omega=(v_i)_{i\ge0}$ in $G$ with $v_0=\rho$, equipped with the product topology.
Since $G$ is locally finite, $\mathcal W_\rho$ is a compact metrizable space:
it is a closed subset of the product of the finite sets $B_G(\rho,i)$.
For $(\omega,\xi)\in\mathcal W_\rho\times\partial\D$, consider the following
condition: there exists $c>0$ such that
\begin{equation}\label{eq:def-B}
h(v_i)\to1,\qquad z(v_i)\to\xi,\qquad
|\xi-z(v_i)|<2e^{-ci}\quad\text{for all sufficiently large $i$.}
\end{equation}
Let $\mathcal E$ be the set of pairs satisfying \eqref{eq:def-B}, and let
$B_*:=\pi_2(\mathcal E)$, where $\pi_2(\omega,\xi)=\xi$.
The set $\mathcal E$ is Borel. Indeed, its first condition can be written as
\[
\bigcap_{k\ge2}\ \bigcup_{M\ge1}\ \bigcap_{i\ge M}
\{(\omega,\xi):h(v_i)>1-1/k\},
\]
and its last condition as
\[
\bigcup_{c\in\Q_+}\ \bigcup_{m\ge1}\ \bigcap_{i\ge m}
\{(\omega,\xi):|\xi-z(v_i)|<2e^{-ci}\}.
\]
The last condition already implies $z(v_i)\to\xi$.
We will use a Borel subset of $B_*$ that has the required exit measure,
rather than make any measurability assertion about the full projection $B_*$.

Let $\nu$ be the law of $((X_i)_{i\ge0},\delta^+)$ on
$\mathcal W_\rho\times\partial\D$ under $\Pp^\rho$.
By L\'evy's $0$--$1$ law, $h(X_i)\to\mathbf 1_A$ almost surely.
Moreover, since consecutive circles intersect, (B) gives, for all large $n$,
\begin{align*}
|\delta^+-z(X_n)|
&\le\sum_{i\ge n}|z(X_{i+1})-z(X_i)|\\
&\le2\sum_{i\ge n}r(X_i)
\le C e^{-c_0n},
\end{align*}
where $C<\infty$ and $c_0>0$ may depend on the realization of the walk.
Choosing $c\in\Q_+$ with $c<c_0$ shows that
\[
\{((X_i)_{i\ge0},\delta^+)\in\mathcal E\}=A
\qquad\text{modulo a $\Pp^\rho$-null set}.
\]
In particular, $\nu(\mathcal E)=\Pp^\rho(A)$.

By inner regularity of the Borel probability measure $\nu$, for each $j\ge1$
choose a compact set $K_j\subset\mathcal E$ such that
$\nu(\mathcal E\setminus K_j)<2^{-j}$, and define
\[
B:=\bigcup_{j\ge1}\pi_2(K_j).
\]
Each $\pi_2(K_j)$ is compact, so $B$ is an $F_\sigma$ subset of $\partial\D$.
Also $B\subset B_*$; in particular, every $\xi\in B$ has a single infinite path
satisfying all the conditions in \eqref{eq:def-B}.

\medskip
\noindent\textbf{Claim 1.} On the event $A$, we have $\delta^+\in B$ almost surely.

\smallskip
\noindent\emph{Proof of Claim 1.}
For every $j$,
\[
\nu\!\left(\mathcal E\setminus\bigcup_{\ell\ge1}K_\ell\right)
\le\nu(\mathcal E\setminus K_j)<2^{-j}.
\]
Thus the measure on the left is zero. Since $A$ agrees almost surely with
$\{((X_i)_{i\ge0},\delta^+)\in\mathcal E\}$, almost every trajectory in $A$
gives a pair in some $K_j$, and hence has $\delta^+\in B$.
\qed
\medskip

\medskip

Fix the counterclockwise orientation on $\partial\D$ and write $I(a,b)$ for the open
counterclockwise arc from $a$ to $b$.
Let
\[
L := I(\delta^-,\delta^+),\qquad
R := I(\delta^+,\delta^-),
\]
and for each $n$ define
\[
p_L(n):=\Pp^{X_n}(\delta^+\in L),
\qquad
p_R(n):=\Pp^{X_n}(\delta^+\in R)=1-p_L(n).
\]

\begin{lemma}\label{lem:balanced}
Almost surely, the events
\[
E_n:=\Big\{\min\big(p_L(n),p_R(n)\big)>\tfrac13\Big\}
\]
occur for infinitely many $n$.
\end{lemma}

\begin{proof}
Let
\[
\mathcal H_n
=\sigma\bigl(G,\Theta,(X_k)_{k\le n}\bigr).
\]
Conditionally on $\mathcal H_n$,
the future exit point $\delta^+$
has the non-atomic exit distribution
$\mu_{X_n}$, while $\delta^-$ is fixed.

Starting from $\delta^-$ in the
counterclockwise direction, define
\[
F_n(\xi)
:=\mu_{X_n}(I(\delta^-,\xi)).
\]
By the probability integral transform,
$F_n(\delta^+)$ is conditionally
uniform on $[0,1]$. Since
$p_L(n)=F_n(\delta^+)$ almost surely,
\[
\mathbb P^{\mathrm{rev}}
(E_n\mid\mathcal H_n)=\frac13.
\]
In particular,
$\mathbb P^{\mathrm{rev}}(E_n)=1/3$.

The sequence $(E_n)_{n\in\mathbb Z}$
is stationary and ergodic under the
time shift of the full marked
bi-infinite environment
$(G,\Theta,(X_k)_{k\in\mathbb Z})$.
Thus $E_n$ occurs infinitely often
almost surely.
\end{proof}

\begin{lemma}\label{lem:hit-future}
On the event $\{\delta^+\in B\}$, almost surely the following holds:
for every $m\ge1$,
\[
\Pp^{X_n}\!\big(\mathrm{hit}\ \{v_m,v_{m+1},\ldots\}\big)>\tfrac14
\quad\text{for infinitely many }n,
\]
where $(v_i)$ is any path witnessing $\delta^+\in B$ as in \eqref{eq:def-B}.
\end{lemma}

\begin{proof}
Fix $m\ge1$ and $n$ with $E_n$.
Let
\[
S_m:=\{\ldots,X_{-2},X_{-1}\}\ \cup\ \{\rho,v_1,\ldots,v_{m-1}\}\ \cup\ \{v_m,v_{m+1},\ldots\}.
\]
Consider the embedded curve $\Gamma_m\subset\overline{\D}$ obtained by joining successive
vertices in the bi-infinite path $(\ldots,X_{-2},X_{-1},\rho,v_1,v_2,\ldots)$ by the geometric
edges of the $\ICP$-carrier (equivalently, by the edges of the fixed plane embedding used
throughout this section). Then $\Gamma_m$ is a connected curve whose only accumulation points
on $\partial\D$ are $\delta^-$ and $\delta^+$, and hence $\Gamma_m$ disconnects $\D$ into two
components whose boundary arcs are $L$ and $R$. Therefore,
\begin{equation}\label{eq:sep-ineq}
\Pp^{X_n}\!\big(\mathrm{hit}\ S_m\big)
\ \ge\ \min\{p_L(n),p_R(n)\}
\ >\ \tfrac13.
\end{equation}

\smallskip

By Lemma~\ref{lem:hit-past-ICP},
\[
\Pp^{X_n}\!\big(\mathrm{hit}\ \{\ldots,X_{-1}\}\big)\ \longrightarrow\ 0.
\]
By Lemma~\ref{lem:hit-finite},
\[
\Pp^{X_n}\!\big(\mathrm{hit}\ \{\rho,v_1,\ldots,v_{m-1}\}\big)\ \longrightarrow\ 0.
\]
Hence along infinitely many $n$ with $E_n$, subtracting these two terms from \eqref{eq:sep-ineq} yields
\[
\Pp^{X_n}\!\big(\mathrm{hit}\ \{v_m,v_{m+1},\ldots\}\big)>\tfrac14
\]
for infinitely many $n$.
\end{proof}

\begin{theorem}\label{thm:Poisson-ICP-detailed}
Under (A) and (B), the unit circle $\partial\D$ realises the Poisson boundary of $G$:
for every invariant event $A$ with $\Pp^\rho(A)>0$, there exists a Borel $B\subset\partial\D$ such that
\[
\Pp^\rho\big(A\ \triangle\ \{\delta^+\in B\}\big)=0.
\]
Equivalently, every bounded harmonic $h$ is of the form
\[
h(v)=\E^v\!\big[g(\delta^+)\big]
\qquad\text{for some bounded Borel }g:\partial\D\to\mathbb{R}.
\]
\end{theorem}

\begin{proof}
Let $A$ and $h(v):=\Pp^v(A)$ be as above,
and let $B$ be the Borel set constructed above.

\smallskip\noindent
\emph{(i) If $(\rho,X_1,\ldots)\in A$, then $\delta^+\in B$ a.s.}
This is Claim~1.

\medskip\noindent
\emph{(ii) If $\delta^+\in B$, then $(\rho,X_1,\ldots)\in A$ a.s.}
Let $(v_i)$ witness $\delta^+\in B$. Choose $m$ so large that
\[
h(v_i)\ \ge\ \tfrac12
\qquad\text{for all } i\ge m.
\]
By Lemma~\ref{lem:hit-future}, there are infinitely many $n$ such that
\[
\Pp^{X_n}\!\big(\mathrm{hit}\ \{v_m,v_{m+1},\ldots\}\big)\ >\ \tfrac14.
\]
Applying \eqref{eq:OS-ineq} with $W=\{v_m,v_{m+1},\ldots\}$ yields, for those $n$,
\[
h(X_n)
\ \ge\ \Pp^{X_n}(\mathrm{hit}\ W)\cdot \inf_{u\in W} h(u)
\ \ge\ \tfrac14\cdot \tfrac12
\ =\ \tfrac18.
\]
Thus $h(X_n)\ge\tfrac18$ infinitely often on $\{\delta^+\in B\}$.
Since by Lévy's $0$--$1$ law $h(X_n)\to \mathbf 1_A$, we conclude $\mathbf 1_A=1$ on $\{\delta^+\in B\}$, i.e.
\[
\Pp^\rho\big(A\ \big|\ \delta^+\in B\big)=1.
\]

\medskip

Combining (i) and (ii) gives $\Pp^\rho\big(A\triangle\{\delta^+\in B\}\big)=0$.
By Proposition~\ref{prop:poisson}, this implies that $\partial\D$ realises the Poisson boundary.
\end{proof}

\medskip

\section{Convergence speed}

\begin{theorem}\label{thm 0.33}
Let $(G, \rho, \Theta)$ be a unimodular $\ICP$-hyperbolic random rooted tame $\IIA$
with $\sup_e\Theta(e)<\pi$, and suppose that $\mathbb{E}[\deg^3(\rho)]<\infty$.
Let $(X_n)$ be a simple random walk on $G$ starting at the root vertex $\rho$, and let
$z(\cdot)$ and $z_h(\cdot)$ denote the Euclidean center and hyperbolic center,
respectively, of the circle corresponding to a vertex in $\D$.
Then, almost surely,
\[
\lim_{n \to \infty} \frac{d_{\mathrm{hyp}}\bigl(z_h(\rho),z_h(X_n)\bigr)}{n}
=
\lim_{n \to \infty} \frac{-\log r(X_n)}{n}
>0.
\]
Moreover, if $(G,\rho,\Theta)$ is ergodic, then the above limit is an almost sure constant.
\end{theorem}

\begin{proof}
We briefly pass from a unimodular random rooted graph to a reversible one by degree-biasing the root.

Let $(G,\rho)$ be a unimodular random rooted graph, and then for every non-negative Borel function
$\varphi$:
\[
  \mathbb{E}\Big[ \sum_{v\in V(G)} \varphi(G,\rho,v) \Big]
  \;=\;
  \mathbb{E}\Big[ \sum_{u\in V(G)} \varphi(G,u,\rho) \Big].
  \tag{MTP}
\]

We now explain how to construct from a unimodular $(G,\rho)$ a new rooted
graph that is reversible for SRW.

Assume that $\mathbb{E}[\deg^3(\rho)]<\infty$.
Define a new probability measure
$\mathbb P^{\mathrm{rev}}$ on rooted
angle-marked IAGs by
\[
\mathbb P^{\mathrm{rev}}(A)
:=
\frac{\mathbb E[
\deg_G(\rho)\mathbf 1
\{(G,\rho,\Theta)\in A\}]}
{\mathbb E[\deg_G(\rho)]}.
\]
We denote expectation with respect to $\mathbb{P}^{\mathrm{rev}}$ by
$\mathbb{E}^{\mathrm{rev}}[\cdot]$.

\begin{lemma}
Under $\mathbb{P}^{\mathrm{rev}}$, the rooted graph $(G,\rho)$ is reversible for
simple random walk.
\end{lemma}

\begin{proof}
Let $(X_n)_{n\ge 0}$ be the simple random walk on $G$ started at $X_0=\rho$.
We must show that for every bounded Borel $F$ on doubly rooted graphs,
$$\mathbb{E}^{\mathrm{rev}}\big[ F(G,X_0,X_1) \big]
  \;=\;
  \mathbb{E}^{\mathrm{rev}}\big[ F(G,X_1,X_0) \big].$$
We first compute $\mathbb{E}^{\mathrm{rev}}[F(G,X_0,X_1)]$ explicitly.
By definition of $\mathbb{P}^{\mathrm{rev}}$,
\[
  \mathbb{E}^{\mathrm{rev}}\big[ F(G,X_0,X_1) \big]
  =
  \frac{1}{\mathbb{E}[\deg(\rho)]}\,
  \mathbb{E}\Big[
    \deg(\rho)\,
    \mathbb{E}\big[ F(G,\rho,X_1) \,\big|\, G,\rho \big]
  \Big].
\]
Conditioned on $(G,\rho)$, $X_1$ is a uniform neighbour of $\rho$, hence
\[
  \mathbb{E}\big[ F(G,\rho,X_1) \,\big|\, G,\rho \big]
  =
  \frac{1}{\deg(\rho)}
  \sum_{v\sim\rho} F(G,\rho,v).
\]
Plugging this into the previous display, the factor $\deg(\rho)$ cancels, giving
\[\mathbb{E}^{\mathrm{rev}}\big[ F(G,X_0,X_1) \big]
  = \frac{1}{\mathbb{E}[\deg(\rho)]}\,
     \mathbb{E}\Big[ \sum_{v\sim\rho} F(G,\rho,v) \Big].
\]
We do the same for the reversed pair $(X_1,X_0)$.
By definition,
\[
  \mathbb{E}^{\mathrm{rev}}\big[ F(G,X_1,X_0) \big]
  =
  \frac{1}{\mathbb{E}[\deg(\rho)]}\,
  \mathbb{E}\Big[
    \deg(\rho)\,
    \mathbb{E}\big[ F(G,X_1,\rho) \,\big|\, G,\rho \big]
  \Big].
\]
Again, conditioned on $(G,\rho)$, $X_1$ is a uniform neighbour of $\rho$, so
\[
  \mathbb{E}\big[ F(G,X_1,\rho) \,\big|\, G,\rho \big]
  =
  \frac{1}{\deg(\rho)}
  \sum_{u\sim\rho} F(G,u,\rho).
\]
Therefore
\[
\mathbb{E}^{\mathrm{rev}}\big[ F(G,X_1,X_0) \big]
  = \frac{1}{\mathbb{E}[\deg(\rho)]}\,
     \mathbb{E}\Big[ \sum_{u\sim\rho} F(G,u,\rho) \Big].
\]
Let
\[
  \varphi(G,u,v) := F(G,u,v)\,\mathbf{1}\{u\sim v\}.
\]
Now apply the mass transport principle to $\varphi$:
\[
  \mathbb{E}\Big[ \sum_{v\in V(G)} \varphi(G,\rho,v) \Big]
  =
  \mathbb{E}\Big[ \sum_{u\in V(G)} \varphi(G,u,\rho) \Big].
\]
Since $\varphi(G,u,v)$ is nonzero only when $u \sim v$, this becomes
\[
  \mathbb{E}\Big[ \sum_{v\sim\rho} F(G,\rho,v) \Big]
  =
  \mathbb{E}\Big[ \sum_{u\sim\rho} F(G,u,\rho) \Big].
\]
Combining the preceding identities proves
the desired reversibility.
\end{proof}

So we may work under the measure where $(G,\rho)$ is reversible for SRW; in particular,
the joint law of $(G,X_n)_{n\in\mathbb{Z}}$ is stationary under time-shifts when the root is degree-biased.
We may also assume $(G,\rho)$ is ergodic; otherwise apply the whole argument on each ergodic component.

By Theorem \ref{thm:icp-poisson}, for our ideal circle pattern $\mathcal{P}$ in $\D$ we have
\[
  z(X_n) \;\to\; \xi^+ \in\partial \D \quad\text{as } n\to+\infty,
\qquad
  z(X_{-n}) \;\to\; \xi^- \in\partial \D \quad\text{as } n\to+\infty,
\]
almost surely, and $\xi^+\neq\xi^-$ almost surely.

Fix one such realization where $\xi^+\neq\xi^-$.
Choose a Möbius transformation
\[
  \Psi : \D \longrightarrow \mathbb{H} := \{w\in\mathbb{C} : \Im w > 0\}
\]
such that
\[
  \Psi(\xi^+)=0,\qquad \Psi(\xi^-)=\infty.
\]
Let $\widetilde{\mathcal{P}} := \Psi(\mathcal{P})$ be the resulting $\ICP$ of $G$ in $\mathbb{H}$.
Denote by $\hat r(v)$ the Euclidean radius of the circle of $v$ in $\widetilde{\mathcal{P}}$,
and by $\hat z(v)\in\mathbb{H}$ its Euclidean centre.
The remaining freedom of Möbius maps preserving $\mathbb{H}$ and fixing $\{0,\infty\}$
is given by real scalings $w\mapsto a w$ with $a>0$.
Such a scaling multiplies all radii $\hat r(v)$ by $a$, but does not change any radius ratio
$\hat r(v)/\hat r(u)$.
Thus, the quantities
\[
  R_n \;:=\; \frac{\hat r(X_n)}{\hat r(X_{n-1})},\qquad n\ge1,
\]
are well-defined and independent of the residual scaling of $\widetilde{\mathcal{P}}$.

Under the reversible, ergodic law, the doubly infinite random walk $(X_n)_{n\in\mathbb{Z}}$
is a stationary ergodic process.
Consequently, the sequence $(\log R_n)_{n\ge1}$ is stationary and ergodic.

Define the flower degree for a vertex $u$,
\[
S(u)\ :=\ \sum_{v\sim u}\deg(v),
\]
the sum of degrees of the neighbours of $u$.

We now check that $\log R_1$ is integrable.
By the refined $\ICP$ ring lemma (Theorem \ref{thm:re-ring-lemma}),
there exists a constant $A=A(\varepsilon_1,\varepsilon_0)>0$ such that for every edge $u\sim v$,
\begin{equation}\label{eq:refined-ring-edge}
\Bigl|\log \frac{\hat r(v)}{\hat r(u)}\Bigr|
\ \le\
A\,\max \{S(u),S(v)\}.
\end{equation}

Consequently,
\[
|\log R_1|
\le A\bigl(S(X_0)+S(X_1)\bigr).
\]
By Lemma~\ref{lem:S-moment} and degree biasing,
\begin{equation}\label{eq:deg-bias}
\mathbb E^{\mathrm{rev}}[S(\rho)]
=
\frac{\mathbb E[\deg_G(\rho)S(\rho)]}
{\mathbb E[\deg_G(\rho)]}
<\infty.
\end{equation}
Stationarity therefore gives
\[
\mathbb E^{\mathrm{rev}}[|\log R_1|]
\le
2A\,\mathbb E^{\mathrm{rev}}[S(\rho)]
<\infty.
\]

Therefore Birkhoff's ergodic theorem applies to $\log R_n$.

Define
\[
  S_n := \sum_{i=1}^n \log R_i
  = \log\hat r(X_n) - \log\hat r(X_0)
  = \log\hat r(X_n) - \log\hat r(\rho).
\]
By Birkhoff's ergodic theorem,
\[
  \frac{1}{n} S_n
  \,=\,
  \frac{1}{n} \sum_{i=1}^n \log R_i
  \,\xrightarrow[n\to\infty]{\text{a.s.}}\,
  \mathbb{E}^{\mathrm{rev}}[\log R_1].
\]
Hence
\begin{equation}\label{eq:upper-radius-rate-refined}
  \lim_{n\to\infty} \frac{-\log\hat r(X_n)}{n}
  =
  -\,\mathbb{E}^{\mathrm{rev}}[\log R_1]
  =:\lambda.
\end{equation}

We first transfer the limit from $\hat r$ back to the radii $r$ in $\D$.
Since $\xi^+\neq \xi^-$, the map $\Psi$ is holomorphic in a neighbourhood of $\xi^+$ in $\widehat{\C}$.
In particular, $\Psi'$ is continuous in a neighbourhood $U$ of $\xi^+$ and bounded above and below on $U$.
On the almost sure event where $z(X_n)\to\xi^+$ and $r(X_n)\to0$, for all sufficiently large $n$
the circle of $X_n$ is contained in $U$.
For such $n$ we have the distortion bounds
\[
r(X_n)\cdot \inf_{w\in B(z(X_n),\,r(X_n))}|\Psi'(w)|
\ \le\ \hat r(X_n)\ \le\
r(X_n)\cdot \sup_{w\in B(z(X_n),\,r(X_n))}|\Psi'(w)|.
\]
By continuity of $\Psi'$ and $r(X_n)\to0$,
\[
\frac{\sup_{w\in B(z(X_n),\,r(X_n))}|\Psi'(w)|}
     {\inf_{w\in B(z(X_n),\,r(X_n))}|\Psi'(w)|}
\ \longrightarrow\ 1,
\qquad
|\Psi'(z(X_n))|\ \longrightarrow\ |\Psi'(\xi^+)|\in(0,\infty).
\]
Hence
\[
\log \hat r(X_n) - \log r(X_n)
=
\log\Bigl(\frac{\hat r(X_n)}{r(X_n)}\Bigr)
\]
is bounded (indeed convergent) and therefore
\[
\lim_{n\to\infty}\frac{-\log r(X_n)}{n}
=
\lim_{n\to\infty}\frac{-\log \hat r(X_n)}{n}
=
\lambda.
\]

By Lemma~\ref{lem:ICP-51},
\[
\limsup_{n\to\infty}
\frac{\log r(X_n)}n<0
\qquad\text{almost surely}.
\]
we conclude that $\lambda>0$.

We now relate $\lambda$ to hyperbolic speed.
For each $n$, since the circle of $X_n$ is contained in $\D$, we have the elementary bound
\begin{equation}\label{eq:r-to-boundary}
r(X_n)\ \le\ 1-|z(X_n)|.
\end{equation}
Also, for any edge $u\sim v$ in $G$, the circles of $u$ and $v$ intersect in $\D$, hence
\[
  |z(u)-z(v)|\;\le\;r(u)+r(v).
\]
Therefore, joining the centres $z(X_n),z(X_{n+1}),\dots$ by straight segments,
\[
  |\xi^+-z(X_n)|
  \le \sum_{i\ge n} |z(X_{i+1})-z(X_i)|
  \le 2\sum_{i\ge n} r(X_i),
\]
and in particular
\begin{equation}\label{eq:dist-upper-refined}
  1-|z(X_n)|
  \le 2\sum_{i\ge n} r(X_i).
\end{equation}
Since $-\frac{1}{n}\log r(X_n)\to\lambda>0$, the tail sum $\sum_{i\ge n} r(X_i)$ has the same exponential order as $r(X_n)$,
and from \eqref{eq:dist-upper-refined} and \eqref{eq:r-to-boundary} we obtain
\[
\lim_{n\to\infty}\frac{-\log(1-|z(X_n)|)}{n}=\lambda.
\]

In the Poincar\'e disc model,
\[
  d_{\mathrm{hyp}}(0,z) = \log\frac{1+|z|}{1-|z|} = -\log(1-|z|)+\log(1+|z|).
\]
Since $|z(X_n)|\to 1$, we have $\log(1+|z(X_n)|)\to\log 2$, hence
\[
  \lim_{n\to\infty}\frac{1}{n}d_{\mathrm{hyp}}(0,z(X_n))
  =
  \lim_{n\to\infty}\frac{-\log(1-|z(X_n)|)}{n}
  =
  \lambda.
\]

It remains to compare the Euclidean and hyperbolic centres.
Let $r_h(v)$ be the hyperbolic radius of the circle of $v$.
We claim that
\[
r_h(v)\le A\,S(v)+C_{\mathrm{ctr}},
\qquad
C_{\mathrm{ctr}}
:=\log\left(1+\frac{2}{\sin^2(\varepsilon_1/2)}\right),
\]
where $A=A(\varepsilon_0,\varepsilon_1)$ is the constant
in the refined Ring Lemma.

To prove this, apply an automorphism of $\D$ sending $z_h(v)$
to $0$. The transformed pattern has the same combinatorics
and intersection angles. The circle of $v$ is now centred at
$0$ with Euclidean radius
\[
s=\tanh\bigl(r_h(v)/2\bigr).
\]
Choose a neighbour $w$ of $v$, and let $t$ be its Euclidean
radius in the transformed pattern and $d$ the distance of
its Euclidean centre from $0$.
The intersection-angle formula gives
\[
d^2=s^2+t^2+2st\cos\Theta(vw).
\]
Since the circle of $w$ is contained in $\D$, we have
$d+t<1$, and consequently
\[
2t\bigl(1+s\cos\Theta(vw)\bigr)<1-s^2.
\]
Set $\kappa_*:=\sin^2(\varepsilon_1/2)>0$.
The bound $\Theta(vw)\le\pi-\varepsilon_1$ implies
\[
1+s\cos\Theta(vw)
\ge \min\{1,1-\cos\varepsilon_1\}
\ge \kappa_*.
\]
Applying the refined Ring Lemma to the transformed pattern,
we also have
\[
t\ge s\exp\bigl(-AS(v)\bigr).
\]
It follows that
\[
\kappa_*\,\sinh r_h(v)
=\kappa_*\,\frac{2s}{1-s^2}
\le \exp\bigl(AS(v)\bigr),
\]
and hence
\[
r_h(v)
\le
\operatorname{arsinh}
\left(\kappa_*^{-1}\exp\bigl(AS(v)\bigr)\right)
\le AS(v)+C_{\mathrm{ctr}}.
\]
This proves the claim.

By Lemma~\ref{lem:S-moment} and degree biasing,
\[
\mathbb E^{\mathrm{rev}}[S(\rho)]
=
\frac{\mathbb E[\deg(\rho)S(\rho)]}
     {\mathbb E[\deg(\rho)]}
\le
\frac{\mathbb E[\deg(\rho)^3]}
     {\mathbb E[\deg(\rho)]}
<\infty.
\]
For every $\eta>0$, stationarity gives
\[
\sum_{n=1}^{\infty}
\mathbb P^{\mathrm{rev}}
\bigl(S(X_n)>\eta n\bigr)
=
\sum_{n=1}^{\infty}
\mathbb P^{\mathrm{rev}}
\bigl(S(\rho)>\eta n\bigr)
<\infty.
\]
Thus the Borel--Cantelli lemma yields
$S(X_n)/n\to0$ almost surely.

The Euclidean centre $z(v)$ belongs to the hyperbolic disk
with centre $z_h(v)$ and radius $r_h(v)$, so
\[
d_{\mathrm{hyp}}(z(v),z_h(v))
\le r_h(v)
\le AS(v)+C_{\mathrm{ctr}}.
\]
Therefore
\(
\frac{d_{\mathrm{hyp}}(z(X_n),z_h(X_n))}{n}
\longrightarrow0
\qquad\text{almost surely}.
\)
The triangle inequality now gives
\[
\left|
\frac{d_{\mathrm{hyp}}(0,z_h(X_n))}{n}
-
\frac{d_{\mathrm{hyp}}(0,z(X_n))}{n}
\right|
\longrightarrow0.
\]
Since the second term converges to $\lambda$, we conclude that
\[
\lim_{n\to\infty}
\frac{d_{\mathrm{hyp}}(0,z_h(X_n))}{n}
=\lambda.
\]
Shifting the basepoint from $0$ to $z_h(\rho)$ changes distances by at most the constant $d_{\mathrm{hyp}}(0,z_h(\rho))$, hence
\[
  \lim_{n\to\infty} \frac{1}{n}d_{\mathrm{hyp}}\bigl(z_h(\rho),z_h(X_n)\bigr)
  = \lambda.
\]
Together with the radius decay rate, we have shown
\[
  \lim_{n\to\infty} \frac{d_{\mathrm{hyp}}\bigl(z_h(\rho),z_h(X_n)\bigr)}{n}
  =
  \lim_{n\to\infty} \frac{-\log r(X_n)}{n}
  =
  \lambda > 0.
\]

If $(G,\rho,\Theta)$ is ergodic, then the limit
\[
  \lambda = \lim_{n\to\infty} \frac{-\log r(X_n)}{n}
\]
is shift-invariant under the stationary bi-infinite marked environment seen from the walk.
By ergodicity it is almost surely constant. Since degree biasing gives a measure
equivalent to the original unimodular law, the almost sure conclusions also
hold under the original law.
\end{proof}

\section{Related topics}\label{sec:further-discussion}
In this final section we list some topics that are closely related to this research.

\paragraph{\emph{\textbf{1. Martin boundary under unbounded degree}}}

For an ergodic unimodular random tame $\IIA$, the bounded-degree condition for identifying the geometric boundary with the Martin boundary seems too strong.  We expect that \(\E[\deg^3(\rho)]<\infty\) is enough. This conjecture is already proposed in \cite{AHNR16} for triangulations.

\begin{conjecture}\label{conj:unbounded}
 Let $(G, \rho, \Theta)$ be an ergodic, unimodular, $\ICP$-hyperbolic random rooted tame $\IIA$ with $\sup_e\Theta(e)<\pi$, and suppose that $\mathbb{E}[\deg^3(\rho)]<\infty$. Then almost surely
 $\partial \D$ is a realization of the Martin boundary of $G$.
\end{conjecture}

\medskip
\paragraph{\emph{\textbf{2. Multi-ended case}}}
In the one-ended setting, the geometry at infinity is connected. For multi-ended unimodular \(\IIA\), new phenomena appear.
First, \(\E[T(\rho)]\) is not expected to determine the \(\ICP\) type on its own. This suggests that one needs additional quantities that record how the geometry is spread across different ends.
Second, random-walk behavior may be more complicated. The walk may pass between ends, which can affect the escape speed and the way the walk approaches the boundary. As a result, the Poisson/Martin/Gromov boundary may be disconnected and can have multiple components.

\medskip
\paragraph{\emph{\textbf{3. Combinatorial Ricci flow}}}

In \cite{GeIIP2}, Ge-Hua-Yu-Zhou proved that if \(\Theta\in(0,\pi)^E\) satisfies \((C_1)\) and there exists a constant
\(c>0\) such that \(T(v)\leq   -c \deg(v)\) for all \(v\in V(G)\), then one can choose an initial metric \(r(0)\) so that the
combinatorial Ricci flow 
\begin{equation}\label{Ricci flow}
\frac{\mathrm{d} r_i}{\mathrm{~d} t}=-K_i \cdot \sinh\left(r_i\right), \quad \forall v_i \in V,
\end{equation}
where $K_i$ is the discrete Gaussian curvature at the vertex $v_i$, converges to a good ideal circle pattern.
In their terminology, a good circle pattern is a circle pattern metric with zero discrete Gaussian curvature.

\begin{conjecture}\label{conj:CRF-positive-drift}
Let \((G,\rho,\Theta)\) be a unimodular ergodic random \(\IIA\) and assume that \(\Theta\in(0,\pi)^E\) satisfies \((C_1)\). Then, the following statements are equivalent:
\begin{enumerate} 
\item[$(1)$] \(\E[T(\rho)]<0\);
\item[$(2)$]  $G$ is almost surely $\ICP$-hyperbolic;
\item[$(3)$] There exists an initial metric \(r(0)\) such that the combinatorial Ricci flow \eqref{Ricci flow} almost surely converges to a good ideal circle pattern.
\end{enumerate}

\end{conjecture}


\section*{Author information}

\begin{multicols}{2}

\AuthorBlock
  {Huabin Ge}
  {Renmin University of China}
  {hbge@ruc.edu.cn}

\AuthorBlock
  {Chuwen Wang}
  {Renmin University of China}
  {chuwenwang@ruc.edu.cn}

\vfill\null
\columnbreak

\AuthorBlock
  {Yangxiang Lu}
  {Renmin University of China}
  {2023000744@ruc.edu.cn}

\AuthorBlock
  {Tian Zhou}
  {Peking University}
  {2201110034@pku.edu.cn}

\end{multicols}


\begin{thebibliography}{99}

\bibitem{AldLy07}
David Aldous and Russell Lyons,
Processes on unimodular random networks,
\emph{Electron. J. Probab.} \textbf{12} (2007), Paper No.~54, 1454--1508.

\bibitem{And70}
E. M. Andreev,
On convex polyhedra of finite volume in Lobachevskii space,
\emph{Math. USSR-Sb.} \textbf{12} (1970), no.~2, 255--259.

\bibitem{ABGN16}
Omer Angel, Martin T. Barlow, Ori Gurel-Gurevich, and Asaf Nachmias,
Boundaries of planar graphs, via circle packings,
\emph{Ann. Probab.} \textbf{44} (2016), no.~3, 1956--1984.

\bibitem{AHNR16}
Omer Angel, Tom Hutchcroft, Asaf Nachmias, and Gourab Ray,
Unimodular hyperbolic triangulations: circle packing and random walk,
\emph{Invent. Math.} \textbf{206} (2016), no.~1, 229--268.

\bibitem{map}
Omer Angel, Tom Hutchcroft, Asaf Nachmias, and Gourab Ray,
Hyperbolic and parabolic unimodular random maps,
\emph{Geom. Funct. Anal.} \textbf{28} (2018), no.~4, 879--942.

\bibitem{AS03}
Omer Angel and Oded Schramm,
Uniform infinite planar triangulations,
\emph{Comm. Math. Phys.} \textbf{241} (2003), no.~2--3, 191--213.

\bibitem{BS96a}
Itai Benjamini and Oded Schramm,
Harmonic functions on planar and almost planar graphs and manifolds via circle packings,
\emph{Invent. Math.} \textbf{126} (1996), no.~3, 565--587.

\bibitem{BS96sq}
Itai Benjamini and Oded Schramm,
Random walks and harmonic functions on infinite planar graphs using square tilings,
\emph{Ann. Probab.} \textbf{24} (1996), no.~3, 1219--1238.

\bibitem{BS01}
Itai Benjamini and Oded Schramm,
Recurrence of distributional limits of finite planar graphs,
\emph{Electron. J. Probab.} \textbf{6} (2001), Paper No.~23, 1--13.

\bibitem{BC13}
Itai Benjamini and Nicolas Curien,
Simple random walk on the uniform infinite planar quadrangulation: subdiffusivity via pioneer points,
\emph{Geom. Funct. Anal.} \textbf{23} (2013), no.~2, 501--531.

\bibitem{IEJ14}
Itai Benjamini, Elliot Paquette, and Joshua Pfeffer,
Anchored expansion, speed and the Poisson--Voronoi tessellation in symmetric spaces,
\emph{Ann. Probab.} \textbf{46} (2018), no.~4, 1917--1956.

\bibitem{BernardiHoldenSun-Percolation-Memoirs2023}
Olivier Bernardi, Nina Holden, and Xin Sun,
Percolation on triangulations: a bijective path to Liouville quantum gravity,
\emph{Mem. Amer. Math. Soc.} \textbf{289} (2023), no.~1440.

\bibitem{BobSpr04}
Alexander I. Bobenko and Boris A. Springborn,
Variational principles for circle patterns and Koebe's theorem,
\emph{Trans. Amer. Math. Soc.} \textbf{356} (2004), no.~2, 659--689.

\bibitem{Buck19}
Ulrike B\"ucking,
On rigidity and convergence of circle patterns,
\emph{Discrete Comput. Geom.} \textbf{61} (2019), no.~2, 380--420.

\bibitem{BBS01}
Dmitri Burago, Yuri Burago, and Sergei Ivanov,
\emph{A course in metric geometry},
Grad. Stud. Math. \textbf{33}, Amer. Math. Soc., Providence, RI, 2001.

\bibitem{Coxeter1950}
H. S. M. Coxeter,
Self-dual configurations and regular graphs,
\emph{Bull. Amer. Math. Soc.} \textbf{56} (1950), 413--455.

\bibitem{N15}
Nicolas Curien,
Planar stochastic hyperbolic triangulations,
\emph{Probab. Theory Related Fields} \textbf{165} (2016), no.~3--4, 509--540.

\bibitem{CLG14}
Nicolas Curien and Jean-Fran\c{c}ois Le Gall,
The Brownian plane,
\emph{J. Theoret. Probab.} \textbf{27} (2014), no.~4, 1249--1291.

\bibitem{CLG16}
Nicolas Curien and Jean-Fran\c{c}ois Le Gall,
The hull process of the Brownian plane,
\emph{Probab. Theory Related Fields} \textbf{166} (2016), no.~1--2, 187--231.

\bibitem{CurMenMir13}
Nicolas Curien, Laurent M\'enard, and Gr\'egory Miermont,
A view from infinity of the uniform infinite planar quadrangulation,
\emph{ALEA Lat. Am. J. Probab. Math. Stat.} \textbf{10} (2013), no.~1, 45--88.

\bibitem{DingGwynne-LQGdim-CMP2020}
Jian Ding and Ewain Gwynne,
The fractal dimension of Liouville quantum gravity: universality, monotonicity, and bounds,
\emph{Comm. Math. Phys.} \textbf{374} (2020), 1877--1934.

\bibitem{GeIIP2}
Huabin Ge, Bobo Hua, Hao Yu, and Puchun Zhou,
Characterization of infinite ideal polyhedra in hyperbolic 3-space via combinatorial Ricci flow,
preprint (2025), arXiv:2506.05036.

\bibitem{GHZ19}
Huabin Ge, Bobo Hua, and Ze Zhou,
Circle patterns on surfaces of finite topological type,
\emph{Amer. J. Math.} \textbf{143} (2021), no.~5, 1397--1430.

\bibitem{GHZ21ICP}
Huabin Ge, Bobo Hua, and Ze Zhou,
Combinatorial Ricci flows for ideal circle patterns,
\emph{Adv. Math.} \textbf{383} (2021), Paper No.~107698, 26 pp.

\bibitem{GeLin2024}
Huabin Ge and Aijin Lin,
The character of Thurston's circle packings,
\emph{Sci. China Math.} \textbf{67} (2024), 1623--1640.

\bibitem{GLWZ2}
Huabin Ge, Yangxiang Lu, Chuwen Wang, and Tian Zhou,
Geometry of unimodular random hyperbolic polyhedra,
preprint (2026), arXiv:2608.03575.

\bibitem{GeIIP1}
Huabin Ge, Hao Yu, and Puchun Zhou,
Infinite ideal polyhedra in hyperbolic 3-space: existence and rigidity,
preprint (2025), arXiv:2506.19528.

\bibitem{Geo16}
Agelos Georgakopoulos,
The boundary of a square tiling of a graph coincides with the Poisson boundary,
\emph{Invent. Math.} \textbf{203} (2016), no.~3, 773--821.

\bibitem{GGN13}
Ori Gurel-Gurevich and Asaf Nachmias,
Recurrence of planar graph limits,
\emph{Ann. of Math. (2)} \textbf{177} (2013), no.~2, 761--781.

\bibitem{GM21}
Ewain Gwynne and Jason Miller,
Random walk on random planar maps: spectral dimension, resistance, and displacement,
\emph{Ann. Probab.} \textbf{49} (2021), no.~3, 1097--1128.

\bibitem{GHS19}
Ewain Gwynne, Nina Holden, and Xin Sun,
Mating of trees for random planar maps and Liouville quantum gravity: a survey,
\emph{Panor. Synth\`eses} \textbf{59} (2023), 41--120.

\bibitem{He96Riem}
Zhengxu He and Oded Schramm,
On the convergence of circle packings to the Riemann map,
\emph{Invent. Math.} \textbf{125} (1996), no.~2, 285--305.

\bibitem{He99}
Zhengxu He,
Rigidity of infinite disk patterns,
\emph{Ann. of Math. (2)} \textbf{149} (1999), no.~1, 1--33.

\bibitem{HS93}
Zhengxu He and Oded Schramm,
Fixed points, Koebe uniformization and circle packings,
\emph{Ann. of Math. (2)} \textbf{137} (1993), no.~2, 369--406.

\bibitem{HeSc95}
Zhengxu He and Oded Schramm,
Hyperbolic and parabolic packings,
\emph{Discrete Comput. Geom.} \textbf{14} (1995), no.~2, 123--149.

\bibitem{HS23}
Nina Holden and Xin Sun,
Convergence of uniform triangulations under the Cardy embedding,
\emph{Acta Math.} \textbf{230} (2023), no.~1, 93--203.

\bibitem{HP17}
Tom Hutchcroft and Yuval Peres,
Boundaries of planar graphs: a unified approach,
\emph{Electron. J. Probab.} \textbf{22} (2017), Paper No.~100, 1--20.

\bibitem{Koebe36}
Paul Koebe,
Kontaktprobleme der konformen Abbildung,
\emph{Ber. Verh. S\"achs. Akad. Wiss. Leipzig Math.-Phys. Kl.} \textbf{88} (1936), 141--164.

\bibitem{LLS25}
Chang Li, Aijin Lin, and Liangming Shen,
The character of ideal circle patterns,
preprint (2026), arXiv:2602.19748.

\bibitem{LyPer16}
Russell Lyons and Yuval Peres,
\emph{Probability on trees and networks},
Cambridge Univ. Press, New York, 2016.

\bibitem{Riv93}
Igor Rivin,
On geometry of convex ideal polyhedra in hyperbolic 3-space,
\emph{Topology} \textbf{32} (1993), no.~1, 87--92.

\bibitem{Riv96}
Igor Rivin,
A characterization of ideal polyhedra in hyperbolic 3-space,
\emph{Ann. of Math. (2)} \textbf{143} (1996), no.~1, 51--70.

\bibitem{RodSul87}
Burt Rodin and Dennis Sullivan,
The convergence of circle packings to the Riemann mapping,
\emph{J. Differential Geom.} \textbf{26} (1987), no.~2, 349--360.

\bibitem{SchlHypCirc05}
Jean-Marc Schlenker,
Hyperideal circle patterns,
\emph{Math. Res. Lett.} \textbf{12} (2005), no.~1, 85--102.

\bibitem{SchlCircSing08}
Jean-Marc Schlenker,
Circle patterns on singular surfaces,
\emph{Discrete Comput. Geom.} \textbf{40} (2008), no.~1, 47--102.

\bibitem{SchlHypPoly}
Jean-Marc Schlenker,
Hyperideal polyhedra in hyperbolic manifolds,
\emph{Math. Res. Lett.} \textbf{20} (2013), no.~4, 773--786.

\bibitem{Schramm1993}
Oded Schramm,
Square tilings with prescribed combinatorics,
\emph{Israel J. Math.} \textbf{84} (1993), no.~1--2, 97--118.

\bibitem{Spring20}
Boris Springborn,
Ideal hyperbolic polyhedra and discrete uniformization,
\emph{Discrete Comput. Geom.} \textbf{64} (2020), no.~1, 63--108.

\bibitem{Steph05}
Kenneth Stephenson,
\emph{Introduction to circle packing: the theory of discrete analytic functions},
Cambridge Univ. Press, Cambridge, 2005.

\bibitem{ThurNotes}
William P. Thurston,
\emph{The geometry and topology of three-manifolds},
Princeton University lecture notes; electronic edition (2002).

\end{thebibliography}
\end{document}